\documentstyle{amltd2004}
\begin{document}
\annalsline{155}{2002}
\received{March 22, 2000}
\startingpage{459}
\def\bye{\end{document}}
 \font\tenrm=cmr10
\input amssym.def
\input amssym.tex
\def\eqref#1{(\ref{#1})}
\def\ritem#1{\item[{\rm #1}]}
\input boxedeps.tex 
\SetepsfEPSFSpecial 
\HideDisplacementBoxes
\def\figin#1#2{\medbreak
$$
 {\BoxedEPSF{#1 scaled
#2}%
}%
$$
\medbreak\noindent}
\catcode`\@=11
\font\twelvemsb=msbm10 scaled 1100
\font\tenmsb=msbm10
\font\ninemsb=msbm10 scaled 800
\newfam\msbfam
\textfont\msbfam=\twelvemsb  \scriptfont\msbfam=\ninemsb
  \scriptscriptfont\msbfam=\ninemsb
\def\msb@{\hexnumber@\msbfam}
\def\Bbb{\relax\ifmmode\let\next\Bbb@\else
 \def\next{\errmessage{Use \string\Bbb\space only in math
mode}}\fi\next}
\def\Bbb@#1{{\Bbb@@{#1}}}
\def\Bbb@@#1{\fam\msbfam#1}
\catcode`\@=12

 \catcode`\@=11
\font\twelveeuf=eufm10 scaled 1100
\font\teneuf=eufm10
\font\nineeuf=eufm7 scaled 1100
\newfam\euffam
\textfont\euffam=\twelveeuf  \scriptfont\euffam=\teneuf
  \scriptscriptfont\euffam=\nineeuf
\def\euf@{\hexnumber@\euffam}
\def\frak{\relax\ifmmode\let\next\frak@\else
 \def\next{\errmessage{Use \string\frak\space only in math
mode}}\fi\next}
\def\frak@#1{{\frak@@{#1}}}
\def\frak@@#1{\fam\euffam#1}
\catcode`\@=12

\newcommand{\R}{{\Bbb R}}
\newcommand{\C}{{\Bbb C}}
\newcommand{\esf}{{\Bbb S}}
\newcommand{\M}{\overline{M}}
\newcommand{\hh}{{\Bbb H}}

\newcommand\escpr[2]{\left<#1,#2\right>}
\newcommand\dbar{\overline{D}}
\newcommand\conn[2]{D\!\raisebox{-2mm}{\it\small #1}\, #2}
\newcommand\rest[2]{\left.#1\right|_{#2}}

 \renewcommand{\leq}{\le}
\renewcommand{\geq}{\ge}

\newenvironment{enum}{\begin{enumerate}
\renewcommand{\labelenumi}{(\roman{enumi})}}{\end{enumerate}}

\let\rr=\R
\let\Sg=\Sigma
\let\sg=\sigma
\let\Ga=\Gamma
\let\th=\theta
\let\ptl=\partial
\let\Om=\Omega
\let\le=\leqslant
\let\ge=\geqslant

 \title{Proof of the Double Bubble Conjecture}


 \acknowledgements{2000 {\it Mathematics Subject Classification}. Primary: 53A10. Secondary: 53C42.}
 \twoauthors{Michael Hutchings, Frank Morgan, Manuel Ritor\'e,}{Antonio Ros}
\shortname{M. Hutchings, F. Morgan, M. Ritor\'e, and A. Ros}
 \institutions{Stanford University, Stanford, CA\\
 {\eightpoint {\it Current addresses\/}:} Institute for Advanced Study, Princeton, NJ and\\
 University of California, Berkeley, CA\\
{\eightpoint {\it E-mail address\/}: hutching@math.ias.edu}\\ \vglue4pt
 Williams College,
Williamstown, MA\\ 
{\eightpoint {\it E-mail address\/}: Frank.Morgan@williams.edu}\\ \vglue4pt
The  University of Granada, Granada, Spain\\
{\eightpoint {\it E-mail address\/}: ritore@ugr.es}
\\ \vglue4pt
The  University of Granada, Granada, Spain\\
{\eightpoint {\it E-mail address\/}:  aros@ugr.es 
\vglue-16pt}}

\centerline{\bf Abstract}
\vglue4pt
We prove that the standard double bubble provides the least-area way to
enclose and separate two regions of prescribed volume in $\rr^3$.
 \vglue-12pt
\section{Introduction}
\vglue-8pt

Archimedes and Zenodorus (see \cite[p.~273]{K}) claimed and Schwarz
\cite{S} proved that the round sphere is the least-perimeter way to
enclose a given volume in $\rr^3$.  The Double Bubble Conjecture, long
believed (see \cite[pp.~300--301]{P}, \cite[p.~120]{B}) but only
recently stated as a conjecture \cite[\S3]{F1}, says that the
familiar double soap bubble of Figure~1, consisting
of two spherical caps separated by a spherical cap or a flat disc,
meeting at 120 degree angles, provides the least-perimeter way to
enclose and separate two given volumes.

\nonumproclaim{Theorem {\rm (see~\ref{th:main})}}
In\/ $\rr^3${\rm ,} the unique perimeter\/{\rm -}\/minimizing double bubble enclosing
and separating regions $R_1$ and $R_2$ of prescribed volumes $v_1$ and
$v_2$ is a standard double bubble as in Figure~{\rm 1,}
consisting of three spherical caps meeting along a common circle at
$120$\/{\rm -}\/degree angles\/{\rm .}\/  $($For equal volumes{\rm ,} the middle cap is a flat
disc.$)$
\endproclaim

The analogous result in $\rr^2$ was proved by the 1990 Williams
College ``SMALL'' undergraduate research Geometry Group \cite{F2}.
The case of equal volumes in $\rr^3$ was proved with the help of a
computer in 1995 by Hass, Hutchings, and Schlafly \cite{HHS},
\cite{Hu}, \cite{HS2} (see \cite{M1}, \cite{HS1},
\cite[Chapt.~13]{M2}).  In this paper we give a complete,
computer-free proof of the Double Bubble Conjecture for arbitrary
volumes in $\rr^3$, using stability arguments, as announced in
\cite{HMRR}.

Reichardt, Heilmann, Lai and Spielman~\cite{RHLS} have generalized our
results to $\rr^4$ and certain higher dimensional cases (when at least
one region is known to be connected).  The 2000 edition of \cite{M2}
treats bubble clusters through these current results.

\figin{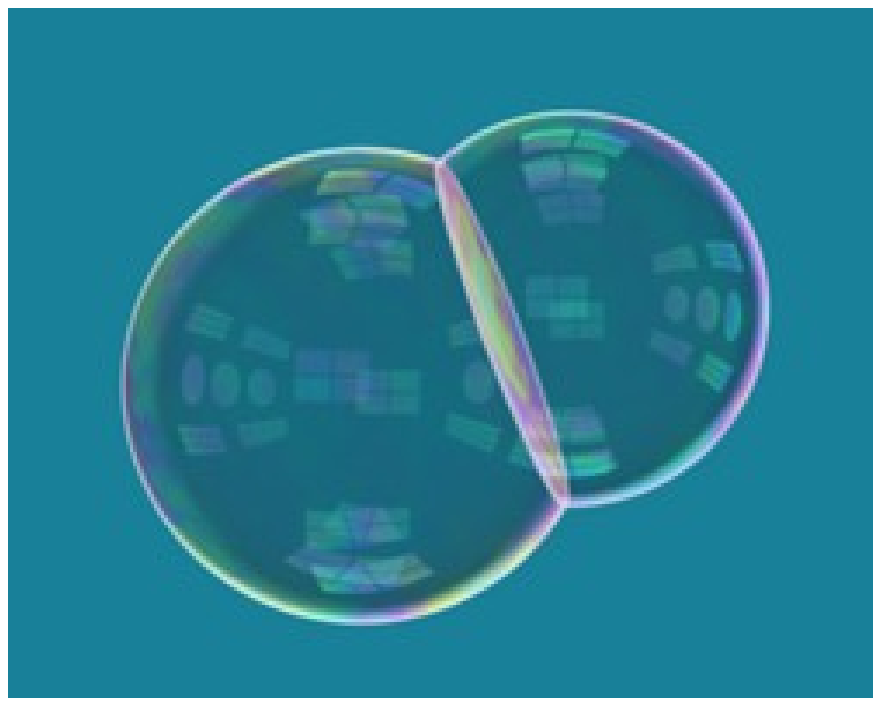}{800}
\begin{quote} Figure 1. The standard double bubble provides the least-perimeter way
to enclose and separate two prescribed volumes.  Computer graphics
copyright John~M. Sullivan, University of Illinois 
(http://\break www.math.uiuc.edu/\~{}jms/Images/) 
\label{fig:standard}
\end{quote}

\demo{Previous results {\rm (see \cite[Chapts.~13 and~14]{M2})}}   Our
strategy for proving Theorem~\ref{th:main} is to assume that a given
double bubble minimizes perimeter and to use this assumption to deduce
that the double bubble is standard.  This strategy is valid only if we
know that a perimeter-minimizing double bubble exists.  F.~Almgren
\cite[Thm.~VI.2]{A} (see \cite[Chapt.~13]{M2}) proved the existence
and almost-everywhere regularity of perimeter-minimizing bubble
clusters enclosing $k$ prescribed volumes in $\rr^{n+1}$, using
geometric measure theory.  J.~Taylor \cite{T} proved that minimizers
in $\rr^3$ consist of smooth constant-mean-curvature surfaces meeting
in threes at 120-degree angles along curves, which in turn could meet
only in fours at isolated points.  An argument suggested by White,
which was written up by Foisy \cite[Thm.~3.4]{F1} and Hutchings
\cite[Thm.~2.6]{Hu}, shows further that any perimeter-minimizing
double bubble in $\rr^{n+1}$ (for $n\ge 2$) has rotational symmetry
about some line.

Unfortunately, the existence proofs depend on allowing the enclosed
regions $R_1$ and $R_2$ to be disconnected.  The complementary
``exterior'' region could also {\it a priori\/} be disconnected.  (If
one tries to require the regions to be connected, they might in
principle disconnect in the minimizing limit, as thin connecting tubes
shrink away.)  Hutchings \cite{Hu} partially dealt with this
complication, using concavity and decomposition arguments to show for
a perimeter-minimizing double bubble that both regions have positive
pressure (see \ref{th:hu}) and hence that the exterior is connected.
Moreover there is a Basic Estimate (see \S \ref{sec:comp}) which
puts upper bounds on the numbers of components of $R_1$ and $R_2$,
depending on the dimension $n$ and the volumes $v_1$, $v_2$.

For equal volumes in $\rr^3$, the Basic Estimate implies that both
enclosed regions are connected.  It can then be shown that a
nonstandard perimeter-minimizing double bubble would have to consist
of two spherical caps with a toroidal band between them
(Fig.~8).  Any such bubble can be described by two
parameters, and Hass and Schlafly \cite{HS2} used a rigorous computer
search of the parameter space to rule out all such possibilities in
the equal volume case, thus proving the Double Bubble Conjecture for
equal volumes in $\rr^3$.  Earlier computer experiments of Hutchings
and Sullivan had suggested that in fact no such nonstandard double
bubbles were stable, and we confirm that in this paper, without using
a computer. \enddemo

\demo{Our proof}  In the present paper we consider arbitrary volumes
$v_1$, $v_2$ in $\rr^3$.  We give a short proof using the Hutchings
Basic Estimate that the larger region is connected
(Proposition~\ref{prop:largeconnected}), and we use a stability
argument (Proposition~\ref{prop:small2comp}) to show that the smaller
region has at most two components, as in Figure~2.
(That the smaller region has at most two components can also be
deduced from the Hutchings Basic Estimate using careful computation;
see \cite[Prop.~4.6]{HLRS}, \cite[14.11--14.13]{M2}.)

\figin{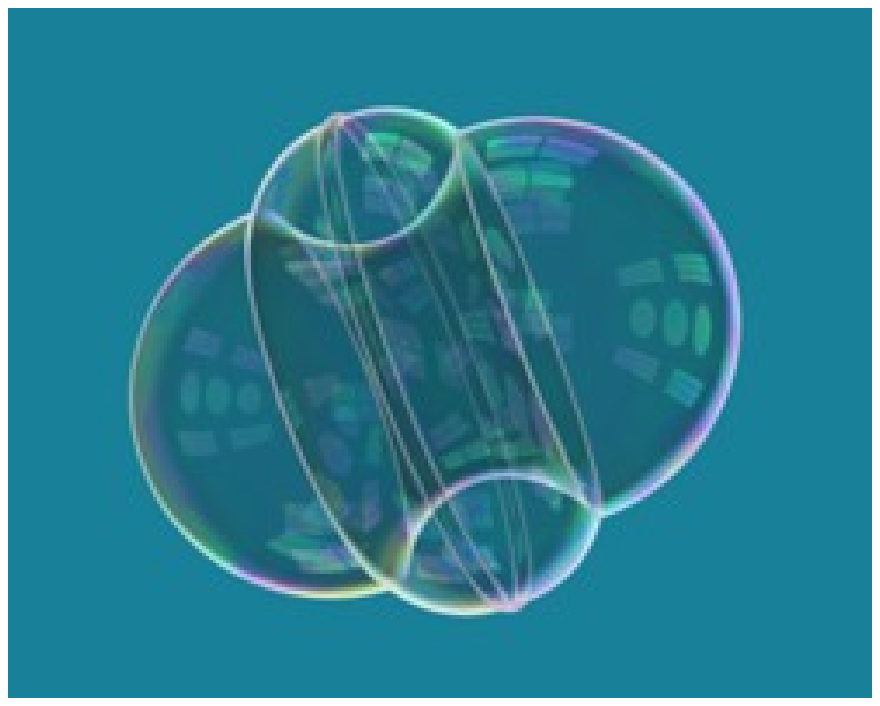}{880}
\begin{quote} Figure 2. A nonstandard double bubble.  One region has two components
(a central bubble and a thin toroidal bubble); the second region is
another toroidal bubble in between.  Computer graphics copyright
John~M. Sullivan, University of Illinois
(http://www.math. uiuc.edu/\~{}jms/Images/) \end{quote}

To prove that an area-minimizing double bubble $\Sg$ is standard,
consider rotations about an axis orthogonal to the axis of symmetry.
At certain places on $\Sg$, the rotation vector field may be tangent to
$\Sg$; i.e., the corresponding normal variation vector field $v$ on
$\Sg$ may vanish.  The axis can be chosen so that these places
separate $\Sg$ into (at least) four pieces
(Proposition~\ref{prop:twocomponents}).  Some nontrivial combinations
$w$ of the restrictions of $v$ to the four pieces vanish on one piece
and respect the two volume constraints.  By stability, $w$ satisfies a
nice differential equation, and hence vanishes on more parts of $\Sg$,
which must therefore be pieces of spheres
(Proposition~\ref{prop:unstable}).  It follows that $\Sg$ must be the
standard double bubble.

The foregoing argument in the proof of Proposition~\ref{prop:unstable}
was inspired by Courant's Nodal Domain Theorem \cite[p.~452]{CH},
which says for example that the first eigenfunction is nonvanishing.
Other applications of this principle to isoperimetric problems and to
the study of volume-preserving stability have been given by Ritor\'e
and Ros \cite{RR}, by Ros and Vergasta \cite{RV}, by Ros and
Souam \cite{RS} and by Pedrosa and Ritor\'e \cite{PR}. \enddemo

{\it Open questions.} We conjecture that the standard double bubble in
$\rr^{n+1}$ is the unique stable double bubble.
Sullivan~\cite[Prob.~2]{SM} has conjectured that the standard
$k$-bubble in $\rr^{n+1}$ $(k\le n+2)$ is the unique minimizer
enclosing $k$ regions of prescribed volume.  This remains open even
for the triple bubble in~$\rr^2$, although Cox, Harrison, Hutchings,
Kim, Light, Mauer and Tilton \cite{CHK} have proved it minimizing in a
category of bubbles with connected regions (which {\it a priori\/} in
principle might bump up against each other).

One can consider the Double Bubble Conjecture in hyperbolic space
$\hh^{n+1}$ or in the round sphere $\esf^{n+1}$.  The symmetry and
concavity results still hold \cite[3.8--3.10]{Hu}.  The case of
$\esf^2$ was proved by Masters~\cite{Ma}.  The cases of $\hh^2$ and
equal volumes in $\hh^3$ and in $\esf^3$ when the exterior is at least
ten percent of $\esf^3$ were proved by Cotton and Freeman~\cite{CF}.

There is also the very physical question in $\rr^3$ of whether the
standard double bubble is the unique stable double bubble with
connected regions.  By our Corollary~\ref{cor:onecomponent}, it would
suffice to prove rotational symmetry.  In $\rr^2$, Morgan and
Wichiramala~\cite{MW} have proved that the standard double bubble is
the unique stable double bubble, except of course for two single
bubbles.

\demo{Contents} Section~2 gives the precise definition of double
bubble and a proof that there is a unique standard double bubble
enclosing two given volumes.  Section~3 provides variational formulas
for our stability arguments.  Section~4 gives some preliminary results
on the geometry of hypersurfaces of revolution with constant mean
curvature (``surfaces of Delaunay'').  Section~5 uses stability
arguments to show that a perimeter-minimizing double bubble in
$\rr^{n+1}$ must be standard if one enclosed region is connected and
the other region has at most two components.  Section~6 proves the
requisite component bounds for perimeter-minimizing double bubbles in
$\rr^3$.  This completes the proof of The Double Bubble Conjecture, as
summarized in Section~7.
\enddemo
\pagebreak

{\it Acknowledgments.} Much of this work was carried out while Morgan
was visiting the University of Granada in the spring of 1999.  Morgan
has partial support from a National Science Foundation grant.
Ritor\'e and Ros have partial support from MCYT research projects BFM2001-3318 and BFM2001-3489, respectively.  We would
like to thank John~M. Sullivan for computer graphics and helpful comments.

\vglue-8pt
\section{Double bubbles}
\vglue-4pt
 
A {\it double bubble\/} in $\rr^{n+1}$ is the union of the topological
boundaries of two disjoint regions of prescribed volumes.  A {\it
smooth double bubble\/} $\Sg\subset\rr^{n+1}$ is a piecewise smooth
oriented hypersurface consisting of three compact pieces $\Sigma_1,
\Sigma_2$ and $\Sigma_0$ (smooth up to the boundary), with a common
$(n-1)$-dimensional smooth boundary $C$ such that $\Sigma_1 +
\Sigma_0$ (resp.\ $\Sigma_2 -\Sigma_0$) encloses a region $R_1$ (resp.\ $R_2$) of prescribed volume $v_1$ (resp.\ $v_2$).  None of these
objects is assumed to be connected.  The unit normal vector field $N$
along $\Sigma$ will be always chosen according to the following
criterion: $N$ points into $R_1$ along $\partial R_1$ and points into
$R_2$ along $\Sigma_2$.  We denote by $\sigma$ and $H$ the second
fundamental form and the mean curvature of $\Sg$.  Note that these
objects are not univalued along the singular set $C$ but they depend
on the sheet $\Sg_i$ we use to compute them.  We will also use the
notation $N_i$, $\sigma_i$ and $H_i$ to indicate the restriction of
$N$, $\sigma$ and $H$ to $\Sigma_i$, $i=0$, $1$, $2$.

Since by Theorem~\ref{th:hu} perimeter-minimizing double bubbles are
smooth double bubbles (geometric measure theory automatically ignores
negligible hair and dirt), throughout the rest of this paper by
``double bubble'' we will mean ``smooth double bubble.''

A {\it standard double bubble} in $\rr^{n+1}$ consists of two exterior
spherical pieces and a separating surface (which is either spherical
or planar) meeting in an equiangular way along a given
$(n-1)$-dimensional sphere $C$.

\proclaim{Proposition}
There is a unique standard double bubble $($up to rigid motions$)$ for
given volumes in $\rr^{n+1}${\rm .}  The mean curvatures satisfy $H_0 =\break H_1 -
H_2${\rm .}
\endproclaim

\demo{Proof}
Consider a unit sphere through the origin and a congruent or smaller
sphere intersecting it at the origin (and elsewhere) at 120 degrees as
in Figure~3.  There is a unique completion to a standard
double bubble.  Varying the size of the smaller sphere yields all
volume ratios precisely once.  Scaling yields all pairs of volumes
precisely once.

The condition on the curvatures follows by plane geometry for $\rr^2$
and hence for $\rr^{n+1}$ (see \cite[Prop.~14.1]{M2}).
\enddemo

\numbereddemo{{R}emark} Montesinos \cite{Mon}  (see \cite[Prob.~2]{SM}) has proved that there
is a unique standard $k$-bubble in $\rr^{n+1}$ for $k\le n+2$.
\enddemo

\begin{center}
\BoxedEPSF{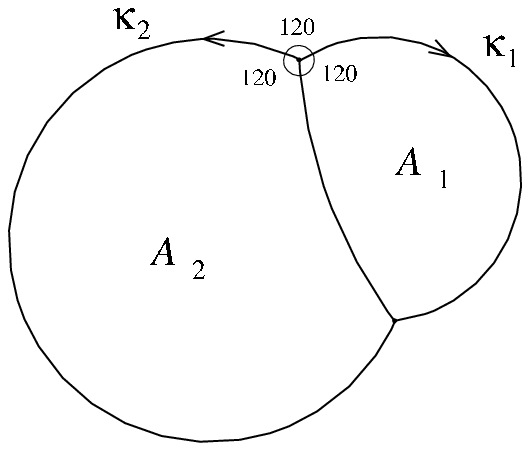 scaled 700} 
\end{center}
\begin{quote} Figure 3. Varying the size of the
smaller cap yields standard double bubbles of all volume ratios.  Then
scaling yields all pairs of volumes.
\end{quote}

\vglue-18pt
\section{Variation formulae}
 
In this section we will consider one-parameter variations
$\{\varphi_t\}_{|t|<\varepsilon} :\Sigma\rightarrow \rr^{n+1}$ of a
double bubble $\Sg\subset\rr^{n+1}$ which are univalued along the
singular set $C$ and when restricted to each one of the pieces
$\Sigma_i$ are smooth (up to the boundary).  Denote by
$X=d\varphi_t/dt$ the associated infinitesimal vector field at $t=0$.
Taking into account our choice of normal vectors to $\Sg$ it is a
standard fact that the derivative of the volume of the regions $R_1$
and $R_2$ are given by
\begin{equation}
\label{eq:1stvolume}
-\int_{\Sg_1}\escpr{X}{N_1}-\int_{\Sg_0}\escpr{X}{N_0}, \quad{\rm and}\quad
-\int_{\Sg_2}\escpr{X}{N_2}+\int_{\Sg_0}\escpr{X}{N_0},\hskip.25in
\end{equation}
respectively.  On the other hand the first derivative of area of the
bubble is given by
$$
\sum_{i=0,1,2}\int_{\Sg_i}\hbox{ div}_{\Sg_i} X,
$$
where $\hbox{ div}_{\Sg_i}$ is the divergence in $\Sg_i$ of a vector
field in $\rr^{n+1}$.  If $\{e_j\}$ is an orthonormal basis of
$T\Sg_i$ and $X$ is a vector field in $\rr^{n+1}$ then
$\hbox{ div}_{\Sg_i} X=\sum_j \escpr{D_{e_j}X}{e_j}$, where $D$ is the
Levi-Civit\'a connection in $\rr^{n+1}$.  As $\hbox{ div}_{\Sg_i}
X=\hbox{ div}_{\Sg_i} X^T-nH_i\escpr{X}{N_i}$, where $X^T$ is the
projection of $X$ to $T\Sg_i$, the Divergence Theorem then implies the
following well-known result.

\proclaimtitle{First variation of area for double bubbles}
\proclaim{Lemma}
\label{lem:1starea}
Let $\Sg\subset\rr^{n+1}$ be a double bubble consisting of smooth
hypersurfaces $\Sg_0${\rm ,} $\Sg_1${\rm ,} $\Sg_2${\rm ,} meeting smoothly along an
$(n-1)$\/{\rm -}\/dimensional submanifold $C${\rm .}  Then the first derivative of the
area along a deformation $\varphi_t(\Sg)$ at $t=0${\rm ,} where $\varphi_t$
is a variation with associated vector field $X${\rm ,} is given by
\begin{equation}
\label{eq:1starea}
-\sum_{i=0,1,2}\int_{\Sg_i} nH_i\escpr{X}{N_i}-\int_C
\escpr{X}{\nu_0+\nu_1+\nu_2},
\end{equation}
where $N_i$ are the normal vectors to the smooth parts $\Sg_i$ of
$\Sg$ and $\nu_i$ are the inner conormals to $C$ inside $\Sg_i${\rm .}
\endproclaim

Suppose that a double bubble $\Sg$ is {\it stationary\/} for any
variation preserving the volume of the regions $R_1$ and $R_2$.  By
Lemma~\ref{lem:1starea} this is equivalent to
\begin{itemize}
\item[(i)] the mean curvatures $H_i$ are constant, with $-H_1+H_2+H_0=0$, and
\item[(ii)] $\nu_0+\nu_1+\nu_2=0$ along $C$.
\end{itemize}

The mean curvature $H_1$ (resp.\ $H_2$) is called the {\it pressure\/}
of the region $R_1$ (resp.\ $R_2$).  From (i) above, we get that if
$H_0>0$, then $R_1$ has larger pressure than $R_2$.

The functions $u_i=\escpr{X}{N_i}$ are the normal components of the
variational field $X$.  If the variation preserves volumes, from
\eqref{eq:1stvolume} they satisfy
\begin{equation}
\label{eq:uvolume}
\int_{\Sg_1} u_1+\int_{\Sg_0} u_0=0, \qquad\qquad \int_{\Sg_2}
u_2-\int_{\Sg_0} u_0=0,
\end{equation}
and, since at the points of the singular set we have $-N_1+N_2+N_0=0$,
we get that
\begin{equation}
\label{eq:uc}
-u_1+u_2+u_0=0 \quad \mbox{\rm along} \, \, C.
\end{equation}

Now we follow the arguments in \cite[Lemma~2.2]{BCE} to show that any
volume preserving infinitesimal variation is integrable.

\proclaim{Lemma}
\label{lem:vp}
Let $\Sg\subset\rr^{n+1}$ be a stationary double bubble{\rm .}  Given smooth
functions $u_i:\Sg_i\to\rr$ such that {\rm \eqref{eq:uvolume}} and
{\rm \eqref{eq:uc}} are satisfied{\rm ,} there is a variation $\{\varphi_t\}$ of
$\Sigma$ which leaves constant the volume of the regions enclosed by
$\varphi_t(\Sg)$ and such that the normal components of the associated
infinitesimal vector field $X$ are the functions $u_i${\rm ,} $i=0,1,2${\rm .}
\endproclaim

\demo{Proof}
The boundary condition \eqref{eq:uc} allows us to construct a smooth
vector field $Z$ on $C$ such that $\escpr{Z}{N_i}=u_i$, which can be
extended smoothly along each $\Sg_i$ so that $\escpr{Z}{N_i}=u_i$.
Let $\{\psi_t\}$ be a one-parameter variation of $\Sg$ associated to
$Z$ (we can take $\psi_t = \psi + t\, Z$).  We choose nonnegative
smooth functions $f_i:\Sg_i\to\rr$, $f_i\neq 0$, with compact support
inside $\hbox{ int}\,\Sg_i$, extended by $0$ to $\Sg$.  For $t$, $s_1$,
$s_2\in \rr$ close to $0$, we consider the three-parameter deformation
$$
\psi_t+s_1\, f_1\,N_1^t+s_2\, f_2\,N_2^t,
$$
where $N_i^t$ is the normal vector to $\psi_t(\Sg_i)$.  Let
$v_i(t,s_1,s_2)$, $i=1$ $2$, be the volume of the deformed region
$R_i$.  Then
$$
\frac{\ptl v_i}{\ptl s_j}(0,0,0)=\int_{\Sg_i} f_j,
$$
which equals $0$ if $i\neq j$ and is positive if $i=j$.  From
conditions \eqref{eq:uvolume}
$$
\frac{\ptl v_i}{\ptl t}(0,0,0)=0.
$$
Applying the Implicit Function Theorem we find smooth functions
$s_1(t)$, $s_2(t)$ with $s_i(0)=0$ such that the volume of the regions
$R_i$ is preserved along the deformation.  Let $X$ be the vector field
associated to this deformation.  Using that the variation is volume
preserving we get that $s_i'(0)=0$.  Hence the normal components of
$X$ are the functions $u_i$.
\enddemo

Now we wish to compute the second derivative of area for a variation
of a double bubble keeping constant the volume of the two enclosed
regions.

\proclaimtitle{Second variation of area for stationary  double bubbles}
\proclaim{Proposition}
\label{prop:2ndarea}
Let $\Sg\subset\rr^{n+1}$ be a stationary double bubble{\rm ,} and let
$\varphi_t$ be a one\/{\rm -}\/parameter variation with associated vector field
$X$ which preserves the volumes of $R_1$ and~$R_2${\rm .}  Then the second
derivative of the area of $\varphi_t(\Sg)$ at $t=0$ is given by
\begin{equation}
\label{eq:2ndarea}
-\int_{\Sg} u\,(\Delta u+|\sg|^2 u)-
\sum_{i=0,1,2}\int_C u_i\,\left\{\frac{\partial
u_i}{\partial\nu_i} + q_iu_i\right\},
\end{equation}
where $u=\escpr{X}{N}${\rm ,} $u_i=\escpr{X}{N_i}${\rm ,} $\Delta$ is the
Laplacian of $\Sg${\rm ,} $|\sg|^2$ is the squared norm of the second
fundamental form{\rm ,} $\nu_i$ is the unit inner normal to $C$ inside
$\Sg_i${\rm ,} and the functions $q_i$ are given by
$q_1=(\kappa_0-\kappa_2)/\sqrt{3}${\rm ,}
$q_2=(-\kappa_1-\kappa_0)/\sqrt{3}$ and
$q_0=(\kappa_1+\kappa_2)/\sqrt{3}$ with $\kappa_i=\sg_i(\nu_i,\nu_i)${\rm ,}
$i=0${\rm ,} $1${\rm ,} $2${\rm .}
\endproclaim

\demo{Proof}
First we recall that the derivative of the mean curvature $H$ is given
by
\begin{equation}
\label{eq:mean}
n\,\frac{dH}{dt} (0) =\Delta u+|\sg|^2u.
\end{equation}

To obtain the second derivative of area we differentiate
\eqref{eq:1starea} with respect to $t$.  The derivative of the
integrals over $\Sg_i$ in \eqref{eq:1starea} equals
\begin{equation}
\label{eq:1stvariation}
-\int_{\Sg} u\,(\Delta u+|\sg|^2u)-
\sum_{i=0,1,2} nH_i\,\left.\frac{d}{dt}\right|_{t=0}\left(\int_{\Sg_i}
\escpr{X}{N_i}\right).
\end{equation}
Let us see that the last sum vanishes.  Let
$a_i=\left.\frac{d}{dt}\right|_{t=0}\left(\int_{\Sg_i}\escpr{X}{N_i}\right)$.
 Since the variation $\varphi_t(\Sg)$ preserves volume, we obtain from
\eqref{eq:1stvolume} that $a_1^{\phantom{|}}+a_0=0$ and $a_2-a_0=0$.  As
$-H_1+H_2+H_0=0$ we conclude
$$
H_1a_1+H_2a_2+H_0a_0=a_0\,(-H_1+H_2+H_0)=0,
$$
which shows that the latter sum in \eqref{eq:1stvariation} vanishes as
we claimed.

It remains to treat the boundary term in \eqref{eq:1starea}.  Since
$\nu_0+\nu_1+\nu_2=0$ on~$C$, differentiating with respect to $t$ we
have
$$
\,\left.\frac{d}{dt}\right|_{t=0}\left(\int_C\escpr{X}{\nu_0+\nu_1+\nu_2}\right)
=\int_C \escpr{X}{\conn{X}{(\nu_0+\nu_1+\nu_2)}}.
$$
Equation \eqref{eq:2ndarea} is then obtained from
Lemma~\ref{lem:2ndarea} below.  To compute $D_{X}\nu_{i}$ the vector
$\nu_{i}$ has been extended as $\nu_{i}^t$ along the integral curves
of $X$, so that $\nu_{i}^t$ is the unit inner conormal to
$\varphi_{t}(C)$ in $\varphi_{t}(\Sg_{i})$.
\enddemo

\numbereddemo{{R}emark}
For bubbles in $\rr^3$, the second variation formula and proof admit
isolated singularities, such as tetrahedral soap film singularities.
For bubbles in $\rr^{n+1}$, $C$ need only be piecewise smooth,
including pieces meeting along an $(n-2)$-dimensional submanifold.  In
addition, the second variation is insensitive to any sets of
${\cal H}^{n-2}$ measure $0$ (see \cite[Lemmas~3.1, 3.3]{MR}).

In a smooth Riemannian ambient manifold $M^{n+1}$, the second
variation has an additional term
$$
- \int_{\Sg} \hbox{ Ric}(N,N)\,u^2
$$
involving the Ricci curvature in the normal direction $N$ (see
\cite[\S 7]{BP}).
\enddemo

\numbereddemo{{R}emark}
By approximation, the second variation formula \eqref{eq:2ndarea}
holds in a distributional sense (see \eqref{eq:defq}) for $X$
piecewise $C^1$ or in $H^1$.
\enddemo

\proclaim{Lemma}
\label{lem:2ndarea}
Under the hypotheses of Proposition {\rm \ref{prop:2ndarea}} we have
\begin{equation}
\label{eq:qi}
 \conn{X}{(\nu_0+\nu_1+\nu_2)}=
\sum_i\left\{\frac{\partial u_i}{\partial\nu_i}+ q_i u_i\right\}\,N_i  =-
\left\{ \frac{d {\th}_2}{d t} (0) N_2 +
\frac{d {\th}_0}{d t} (0) N_0\right\},
\end{equation}
where $\theta_2=\angle (\nu_1,\nu_2)${\rm ,} $\th_0=\angle(\nu_1,\nu_0)$ are
the angles determined by the sheets of $\varphi_t(\Sg)$ along the
singular set{\rm .}
\endproclaim

\demo{Proof}
Let $Y = X^C$ be the orthogonal projection of $X$ to the tangent
bundle $TC$.  For each $i$ we have
$X=Y+\escpr{X}{\nu_i}\,\nu_i+\escpr{X}{N_i}\,N_i$.  As $Y$ is tangent
to $C$ we have $\escpr{Y}{D_X(\sum_i\nu_i)}=0$.  We also have
$\escpr{\nu_i}{D_X\nu_i}=0$ for each $i$ since $|\nu_i|=1$.  Moreover
$$
\escpr{D_X \nu_i}{N_i}=-\escpr{\nu_i}{D_X
N_i}=\sg_i(\nu_i,X^i)+\frac{\partial u_i}{\partial\nu_i},
$$
where $X^i$ is the orthogonal projection of $X$ to $T\Sg_i$.  Hence
$$
D_X (\nu_0+\nu_1+\nu_2)=\sum_i\left\{\frac{\partial
u_i}{\partial\nu_i}+\sg_i(\nu_i,\nu_i)\escpr{X}{\nu_i}+\sg_i(Y,\nu_i)
\right\}\,N_i.
$$

Observe that $\escpr{D_{Y}\nu_i}{\nu_i}=0$, and that $\sum_i
D_{Y}\nu_i=D_{Y}(\sum_{i}\nu_{i})$ vanishes since $Y$ is tangent to
the singular set, where $\sum_{i}\nu_{i}$ is identically zero.  Of
course this implies
$\sum_{i}(D_{Y}\nu_{i})^{C}=(D_{Y}(\sum_{i}\nu_{i}))^{C}=0$.  Hence we
see that
\begin{eqnarray*}
0 = \sum_i D_{Y}\nu_i &=&\sum_i\left\{(D_{Y}\nu_i)^C
+\escpr{D_{Y}\nu_i}{\nu_i}\,\nu_i
+\escpr{D_{Y}\nu_i}{N_i}\,N_i
\right\} \\
&=& \sum_i \escpr{D_{Y}\nu_i}{N_i}\,N_i = \sum_i\sg_i(Y,\nu_i)\,N_i .
\end{eqnarray*}
Therefore we get that $D_X (\nu_0+\nu_1+\nu_2)=\sum_i\left\{\frac{\partial
u_i}{\partial\nu_i}+\kappa_i \escpr{X}{\nu_i}
\right\}\,N_i.$

Taking into account that
\begin{equation}
\nu_1=\frac{1}{\sqrt{3}}\,(N_0-N_2),\qquad
\nu_2=\frac{1}{\sqrt{3}}\,(-N_1-N_0),\qquad
\nu_0=\frac{1}{\sqrt{3}}\,(N_1+N_2),
\label{eq:nui}
\end{equation}
we obtain that
$$
\escpr{X}{\nu_1}=\frac{1}{\sqrt{3}}\,(u_0-u_2),\quad
\escpr{X}{\nu_2}=\frac{1}{\sqrt{3}}\,(-u_1-u_0),\quad
\escpr{X}{\nu_0}=\frac{1}{\sqrt{3}}\,(u_1+u_2).
$$
As $-N_1+N_2+N_0=0$ we have
\begin{eqnarray*}
\sqrt{3}\,\sum_i \kappa_i\escpr{X}{\nu_i}\,N_i=\sqrt{3}\sum_i q_i u_i N_i
\ +&&\hskip-17pt\left\{
(\kappa_0 u_2-\kappa_2 u_0)\,N_1
+(\kappa_1 u_0-\kappa_0 u_1)\,N_2\right.
\\
&&\hskip-7pt \left.+(\kappa_2 u_1-\kappa_1 u_2)\,N_0
\right\}.
\end{eqnarray*} 
The summand between brackets is a vector whose coordinates coincide,
up to sign, with the determinant of the matrix
$$
\left( \begin{array}{ccc}
N_1^i & \kappa_1 & u_1 \\
N_2^i & \kappa_2 & u_2 \\
N_0^i & \kappa_0 & u_0
\end{array}\right),
$$
where $N_j=(N_j^1,\ldots,N_j^n)$.  But this determinant vanishes since
$-N_1^i+N_2^i+N_0^i=0$, $-u_1+u_2+u_0=0$, and
$-\kappa_1+\kappa_2+\kappa_0=0$.  This last equality holds because
$-H_1+H_2+H_0=0$, and $-\sg_1 (Z,T)+\sg_2 (Z,T)+\sg_0 (Z,T)=0$ for any
vector $Z$ and $T$ in $TC$.  Hence the first part of \eqref{eq:qi}
follows.

To prove the remaining part we write $\nu_2=R(\th_2)\nu_1$,
$\nu_0=R(\th_0)\nu_1$, where $R(\th )$ is the rotation in the plane
spanned by $\nu_1$ and $N_1$ given, in this basis, by
$$
R(\th)=
\left( \begin{array}{ccc}
\cos\th & -\sin\th \\
\sin\th & \cos\th
\end{array}\right).
$$
We have
$$
D_X (\nu_0+\nu_1+\nu_2)=(\hbox{ Id}+R(\th_2)+R(\th_0) )\frac{d {\nu}_1}{d t}
+ \frac{d{\th}_2}{d t} (- N_2)+\frac{d {\th}_0}{d t} (- N_0),
$$
and \eqref{eq:qi} follows since the first summand vanishes at $t =0$.
\enddemo

\numbereddemo{{R}emark}
For a variation such that the angles of the sheets are preserved, we
have $D_X (\nu_0+\nu_1+\nu_2)=0$ (since $\nu_0+\nu_1+\nu_2=0$ for all
$t$), so by \eqref{eq:qi}, the boundary term in the second variation
formula \eqref{eq:2ndarea} vanishes.
\enddemo

Consider a stationary bubble $\Sg$.  We say that a function
$u:\bigcup\Sigma_i\rightarrow\R$ defined on the disjoint union of the
$\Sigma_i$ is {\it admissible\/} if the restrictions $u_i$ to the
smooth pieces $\Sg_i$ of $\Sg$, lie in the Sobolev space $H^1$ (of
functions in $L^2$ whose gradient is squared integrable) and satisfy
the boundary condition
\begin{equation}
\label{boundcond}
u_1 = u_2 + u_0 \quad \hbox{ along the singular set } C.
\end{equation}

The space of admissible functions satisfying the zero mean value
conditions
\begin{equation}
\label{mvcond}
\int_{\Sg_1} u_1+\int_{\Sg_0}u_0 = \int_{\Sg_2} u_2-\int_{\Sg_0} u_0 = 0
\end{equation}
will be denoted by ${\cal F} (\Sigma)$.  From the results at the
beginning of this section, we see that admissible functions correspond
to deformations of $\Sigma$ and that ${\cal F} (\Sigma)$ are the
infinitesimal variations of those deformations which preserve the
volume of the regions $R_1$ and $R_2$.  The bilinear form on the space
of admissible functions for the second variation of the area
\eqref{eq:2ndarea} will be denoted by $Q$, and it is given by
\begin{eqnarray}
\label{eq:defq}
\qquad Q(u,v) &=& \int_{\Sigma} \left\{\escpr{\nabla u}{\nabla v}-|\sg|^2 u v
\right\} -
\sum_{i=0,1,2}\int_{C} q_i u_i v_i
\\
 &=& - \int_{\Sigma} ( \Delta u + |\sg|^2 u )\,v \, -\sum_{i=0,1,2}
 \int_{C} \left\{ \left(\frac{\partial u_i}{\partial\nu_i} +q_i
u_i\right)\,v_i \right\},\nonumber
\end{eqnarray}
where $\nu_i$ is the inner normal to $C$ inside $\Sg_i$ and $q_i$ are
the functions defined in the statement of Proposition
\ref{prop:2ndarea}.  We will say that a (smooth) double bubble
$\Sigma$ is {\it stable\/} if it is stationary and $Q(u,u)\ge 0$ for
any $u\in {\cal F} (\Sigma)$.  We shall say that it is {\it
unstable\/} if it is not stable.  By Lemma~\ref{lem:vp} a
perimeter-minimizing double bubble is stable.

\proclaim{Lemma}
\label{Qe0}
Let $\Sigma$ be a stable double bubble and $u\in {\cal F} (\Sigma)$
such that $Q(u,u) = 0${\rm .}  Then $u$ is smooth on the interior of
$\Sigma_i${\rm ,} $i=0${\rm ,} $1${\rm ,} $2${\rm ,} and there exist real numbers $\lambda_0${\rm ,}
$\lambda_1${\rm ,} $\lambda_2${\rm ,} with $\lambda_1 = \lambda_0 + \lambda_2${\rm ,}
such that
$$
\Delta u_i + |\sigma|^2 u_i   = \lambda_i, \  \hbox{ on}\ \Sigma_i.
$$
\endproclaim

\demo{Proof}
The stability of $\Sigma$ implies that $Q(u+tv,u+tv)\ge 0$ for any
$v\in {\cal F}$ and $t\in\rr$.  Therefore $Q(u,v)=0$ and so, taking
arbitrary functions with mean zero and support inside the interior of
$\Sg_i$ we conclude that the displayed equation holds in a
distributional sense.  From elliptic regularity, $u$ is smooth on the
interior of $\Sg_i$.
\enddemo

A smooth admissible function $u$ is said to be a {\it Jacobi function}
if it corresponds to an infinitesimal deformation of $\Sigma$ which
preserves the mean curvature of the pieces $\Sigma_k$, and the fact
that these pieces meet in an equiangular way along its singular set.
By formulae \eqref{eq:mean} and \eqref{eq:qi}, we have that $u$ is a
Jacobi function if and only if
\begin{equation}
\label{Jacobi}
\left\{ \begin{array}{ll}
\Delta u + |\sg|^2 u = 0, &\hbox{ on}\ \Sigma,
\\
\displaystyle-\left(\frac{\partial u_1}{\partial\nu_1} +q_1 u_1\right) =
\frac{\partial
u_2}{\partial\nu_2} + q_2 u_2 = \frac{\partial u_0}{\partial\nu_0} + q_0
u_0, & \hbox{ along}\ C.
\end{array}\right.
\end{equation}
Any Killing vector field $Y$ of $\R^{n+1}$ gives a Jacobi function on
$\Sigma$, $u= \langle Y, N\rangle$. \pagebreak

\proclaim{Lemma}
\label{Qww}
Let $S\subset\Sigma$ be a subdomain with piecewise\/{\rm -}\/smooth boundary and
$u$ a Jacobi function on $\Sigma$ which vanishes on $\partial S$ $($in
particular we assume that all the $u_i$ vanish at $\ptl S\cap C).$
If $w$ is defined by
$$
w =\left\{ \begin{array}{ll}
u, & \hbox{ on}\ S, \\
0, & \hbox{ on}\ \Sigma - S,
\end{array}\right.
$$
then $w$ is an admissible function and $Q(w,w) = 0$\/{\rm .}\ 

Let $S'\subset\Sigma$ be a second subdomain{\rm ,} with the same properties
of $S${\rm ,} and $w'$ its associated admissible function{\rm .}  If the interiors
of $S$ and $S'$ are disjoint{\rm ,} then $Q(w,w^\prime)=0${\rm .}
\endproclaim

{\it Proof}.
A Jacobi function $u$ satisfies $Q(u,v)=0$ for any admissible function
$v$ by \eqref{eq:defq} and \eqref{Jacobi}.  Therefore the equalities
$Q(w,w) = Q(u,w) = 0$ prove the first assertion.  The second one is
trivial.
\hfill\qed

\section{Area-minimizing double bubbles and Delaunay hypersurfaces}
 
As described in the Previous Results section of the introduction,
F.~Almgren \cite[Thm.~VI.2]{A} (see \cite[Chapt.~13]{M2}) proved the
existence and almost-everywhere regularity of perimeter-minimizing
bubble clusters enclosing $k$ prescribed volumes in $\rr^{n+1}$, using
geometric measure theory.  Using symmetry, concavity, and
decomposition arguments, Hutchings analyzed the structure of
minimizing double bubbles.

\proclaim{Theorem}
\label{th:hu}
\begin{itemize}
\item[{\rm (a)}] {\rm (after White, \cite[Thm.~3.4]{F1}, \cite[Thm.~2.6]{Hu}}$).$
An area\/{\rm -}\/minimizing double bubble in $\rr^{n+1}$ $($for $n\ge 2)$ is a
hypersurface of revolution about some line $L${\rm .}
\item[{\rm (b)}] {\rm (\cite[Cor.~3.3]{Hu}).} In an area\/{\rm -}\/minimizing double bubble{\rm ,}
both enclosed regions have positive pressure{\rm .}
\item[{\rm (c)}] {\rm (\cite[Thm.~5.1]{Hu}).} An area\/{\rm -}\/minimizing double bubble is
either the standard double bubble or consists of a topological sphere
with a finite tree of annular bands attached as in
Figure~{\rm 4.} The two caps are pieces of spheres{\rm ,} and the root
of the tree has just one branch{\rm .}  All pieces are smooth $($Delaunay$)$
hypersurfaces meeting in threes at $120$\/{\rm -}\/degree angles along\break round
$(n-1)$\/{\rm -}\/spheres{\rm .}
\end{itemize}

\endproclaim

\phantom{hi}
\vglue-24pt
Let $\Om$ a connected component of the regions $R_1$ or $R_2$ in a
nonstandard minimizing double bubble.  Then either the smooth pieces
in the boundary of $\Om$ are all annuli or $\ptl\Om$ is the union of
two spherical caps $D_1$ and $D_2$ and one annulus $M_0$.  In the
first case we shall refer to $\Om$ as a {\it torus component}, and in
the latter one as the {\it spherical component}.

\begin{center}\BoxedEPSF{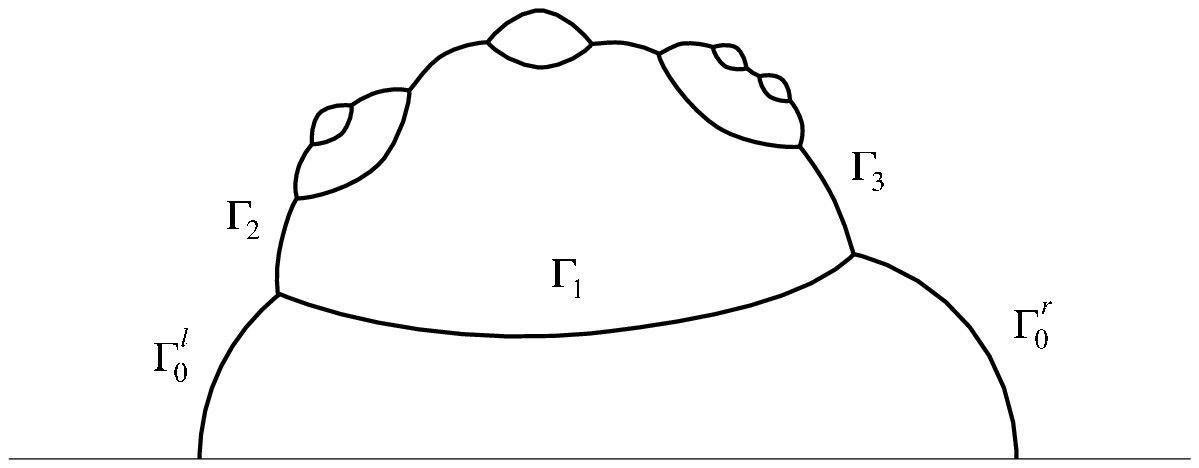 scaled 700} \end{center}
\begin{quote} Figure 4. Generating curve of nonstandard area-minimizing double bubble in
$\rr^{n+1}$.
\end{quote}
 
Now we review some facts about hypersurfaces of revolution with
constant mean curvature in $\rr^{n+1}$, known as Delaunay
hypersurfaces (see \cite{D}, \cite{E}, \cite{Ke} on $\rr^3$ and
\cite{Hs} on $\rr^{n+1}$).  Let $\Sg\subset\rr^{n+1}$ be a
hypersurface and assume that $\Sg$ is invariant under the action of
the group $O(n)$ of isometries of $\rr^{n+1}$ fixing the $x_1$-axis.
The hypersurface $\Sg$ is generated by a curve $\Ga$ contained in the
$x_1x_2$-plane.  The coordinates $x_1$, $x_2$, will be denoted by $x$,
$y$, respectively.  We parametrize the curve $\Ga=(x,y)$ by
arc-length $s$.  If $\alpha$ is the angle between the tangent to $\Ga$
and the positive $x$-direction we shall always choose the normal
vector field $N=(\sin\alpha, -\cos\alpha)$.  Then we have

\proclaim{Lemma}
\label{lem:eqns}
The generating curve $\Ga$ of an $O(n)$\/{\rm -}\/invariant hypersurface
$\Sg\subset\rr^{n+1}$ with mean curvature $H$ with respect to the
normal vector $N=(\sin\alpha,-\cos\alpha)$ satisfies the following
system of ordinary differential equations
\begin{eqnarray}
 x' &=& \cos\alpha, \label{eq:generating}
\\
y' &= &\sin\alpha,
\nonumber\\
 \alpha' &= &-nH + (n-1)\,\frac{\cos\alpha}{y}.\nonumber
\end{eqnarray}
Moreover{\rm ,} if $H$ is constant then the above system has the first
integral
\begin{equation}
\label{eq:fint}
y^{n-1}\cos\alpha- Hy^n=F.
\end{equation}
\endproclaim

The constant $F$ in \eqref{eq:fint} is called the {\em force\/} of the
curve $\Ga$.  Existence of the first integral is standard in the
Calculus of Variations (see \cite[\S 3.4]{GH} and the references
therein).  For constant mean curvature surfaces see
\cite[pp.~138--139]{P}, with earlier reference to Beer and \cite[\S
3]{KKS}.

From Lemma~\ref{lem:eqns} we can obtain the following known properties.

\proclaim{Proposition}
\label{prop:solutions}
Any local solution of the system {\rm \eqref{eq:generating}} is a part of a
complete solution $\Ga${\rm ,} which generates a hypersurface $\Sg$ with
constant mean curvature of several possible types $($see
Figure~{\rm 5}$)$. \pagebreak
\begin{itemize}
\ritem{(i)} If $FH>0$ then $\Ga$ is a periodic graph over the $x$-axis{\rm .} It
generates a periodic embedded unduloid{\rm ,} or a cylinder{\rm .}

\ritem{(ii)} If $FH<0$ then $\Ga$ is a locally convex curve and $\Sg$ is a
nodoid{\rm ,} which has self\/{\rm -}\/intersections{\rm .}

\ritem{(iii)} If $F=0$ and $H\neq 0$ then $\Sg$ is a sphere{\rm .}

\ritem{(iv)} If $H=0$ and $F\neq 0$ we obtain a catenary which
generates an embedded catenoid $\Sg$ with $F>0$ if the normal
points down and $F<0$ if the
normal points up{\rm .}

\ritem{(v)} If $H=0$ and $F=0$ then $\Ga$ is a straight line orthogonal to the
$x$\/{\rm -}\/axis which generates a hyperplane{\rm .}

\ritem{(vi)} If $\Sigma$ touches the $x$\/{\rm -}\/axis{\rm ,} then it must be a sphere or a
hyperplane{\rm .}

\ritem{(vii)} The curve $\Ga$ is determined{\rm ,} up to translation along the
$x$\/{\rm -}\/axis{\rm ,} by  the pair $(H,F)${\rm .}
\end{itemize}

\endproclaim

The generating curves of nodoids and unduloids are called {\it
nodaries\/} and {\it undularies}.  Since we shall often identify the
curves and the generated hypersurfaces we shall refer to them as
nodoids and unduloids.
\figin{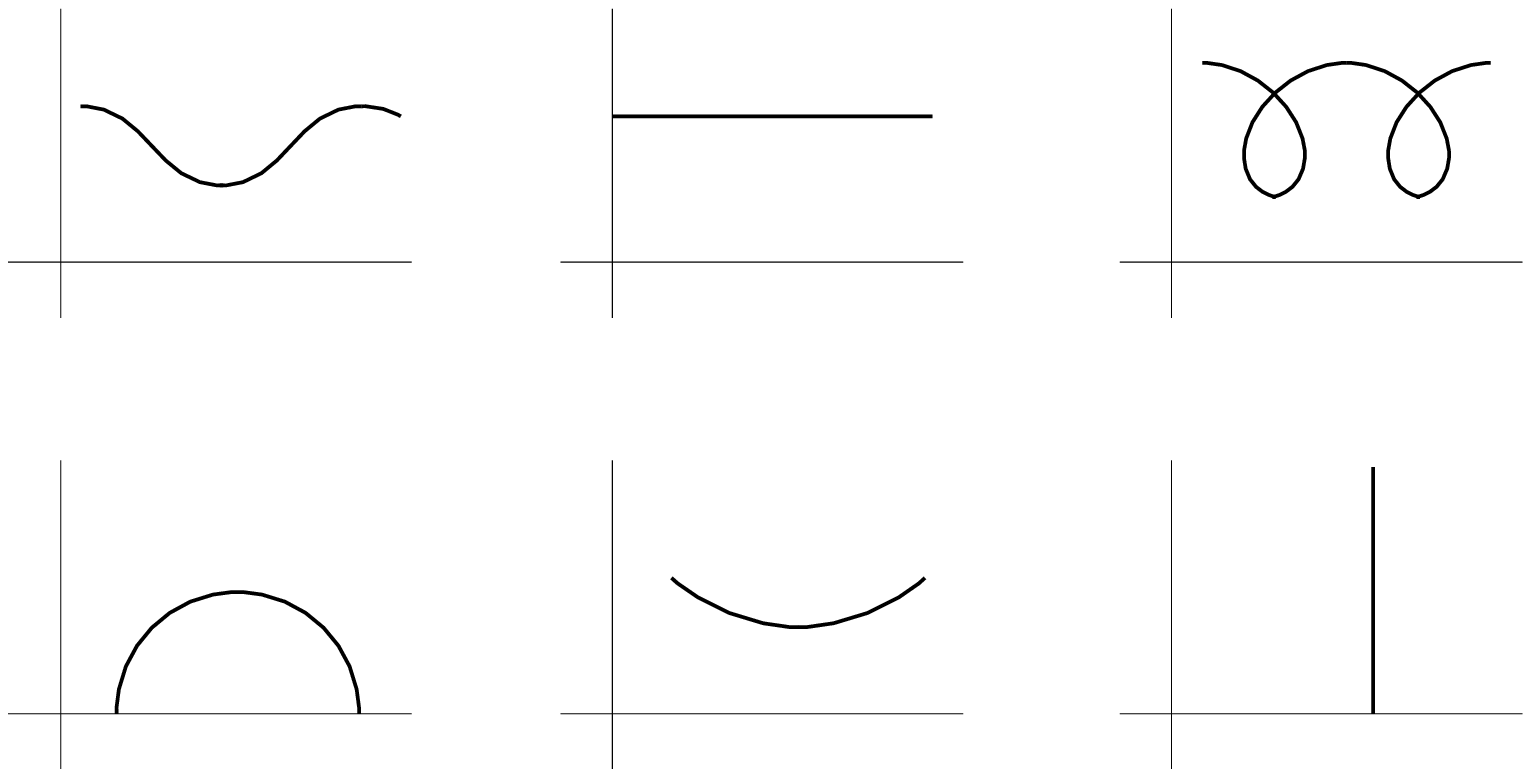}{550}
\begin{quote} Figure 5. The generating curves for Delaunay hypersurfaces: unduloid,
cylinder, nodoid, sphere, catenoid, hyperplane.
\end{quote}
 \vglue4pt

\numbereddemo{{R}emark}
Henceforth we shall use the following properties of generating curves
of Delaunay hypersurfaces.
\begin{itemize}
\item[(i)] Unduloids and nodoids have positive mean curvature with respect
to the normal which points downward at the maximum of the
$y$-coordinate.
\item[(ii)] The nodoid is convex in the sense that the normal vector rotates
monotonically.  This follows from equations \eqref{eq:generating} and
\eqref{eq:fint}.
\end{itemize}

\enddemo

\proclaimtitle{Force balancing {\cite[(3.9)]{KKS}}}
\proclaim{Lemma}
\label{lem:forcebalancing}
Assume that three generating curves $\Ga_i${\rm ,} $i=0${\rm ,} $1${\rm ,} $2$ of
hypersurfaces with constant mean curvature $H_i$ and forces $F_i$ meet
at some point{\rm .}  If $-H_1+H_2+H_0=0$ and $-N_1+N_2+N_0=0$ at this
point{\rm ,} then
\begin{equation}
\label{eq:energy}
-F_1+F_2+F_0=0.
\end{equation}
\endproclaim

The lemma follows directly from \eqref{eq:fint}.

\proclaim{Lemma}
\label{lem:forcegamma1}
Let $\Sg$ be a nonstandard minimizing double bubble in $\rr^{n+1}${\rm ,} as
in Figure~{\rm 4,} and let $R_1$ be the region of larger or equal
pressure{\rm .}  Assume that the spherical component $\Om$ is contained in
$R_1${\rm .}  Let $\Ga_1$ be the generating curve of $M_0=\Sg_0\cap\ptl\Om${\rm .}

Then the force of $\Ga_1$ is positive and $\Ga_1$ is an unduloid or
catenoid and in particular a graph{\rm .}
\endproclaim

{\it Proof}.
Let $\Ga_0^l$, $\Ga_0^r$ be the left and right circles in $\ptl\Om$.
Consider the embedded curve $\Ga$ determined by $\Ga_0^r$, $\Ga_0^l$
and $\Ga_1$.  Let $\Ga_2$ be the third generating curve meeting
$\Ga_0^l\cap\Ga_1$ and $\Ga_3$ the one meeting $\Ga_0^r\cap\Ga_1$.

If the force of $\Ga_1$ is negative then $\Ga_1$ is a nodoid with
positive mean curvature since $H_0=H_1-H_2\ge 0$.  The graph $\Ga$ is
convex and, since $\Ga$ meets $L$ orthogonally, its total curvature
equals $\pi$.  At each one of the vertices $\Ga_0^r\cap\Ga_1$,
$\Ga_0^l\cap\Ga_1$ the inner angle of $\Ga$ is exactly $\pi/3$.  By
force balancing~\ref{lem:forcebalancing}, both $\Ga_2$ and $\Ga_3$
have positive force and they are unduloids.  Since $R_2$ has positive
pressure both $\Ga_2$ and $\Ga_3$ are inward graphs with respect to
$\Ga_1$ (i.e., the exterior region lies above $\Gamma_1$ and above
$\Gamma_2$).  Hence the two circular arcs $\Ga_0^l$, $\Ga_0^r$ have
angular measure larger than $\pi/3$.  So the total curvature of
$\Gamma$ is larger than $4\pi/3$, which is a contradiction.

If the force of $\Ga_1$ is $0$ then $\Ga_1$ is part of a circle or of
a line orthogonal to the axis of revolution $L$.  The former
possibility is discarded by the same argument used for nodoids.  The
latter possibility is clearly not possible.

Hence the force of $\Ga_1$ is positive and $\Ga_1$ is an unduloid or
catenoid and graph.
\hfill\qed\vglue4pt

By a similar argument to the one used in Lemma~\ref{lem:forcegamma1}
we obtain

\proclaim{Lemma}
\label{lem:gammas}
Let $\Sg$ be a double bubble of revolution such that both regions
have positive pressure{\rm .} Then it is not possible that $\Sg$ contains
pieces of spheres $\Ga_0^l${\rm ,} $\Ga_0^r${\rm ,} $\Ga_1${\rm ,} $\Ga_2${\rm ,} $\Ga_3$ as
in Figure~{\rm 6,} with $\Ga_1\subset\Sg_0${\rm .}
\endproclaim

\proclaim{Lemma}
\label{lem:delaunaygraph}
Let $\Sg$ be a nonstandard minimizing double bubble in $\rr^3${\rm .}  Let
$\th_i$ be the subtending angle of the spherical caps $D_i$ as in
Figure~{\rm 7.}
 \vglue1pt
 {\rm (i)} \hangindent=39pt\hangafter=1
 If $\th_1$, $\th_2\le\pi/6$ then $\th_1=\th_2$ and $M_0$ is
symmetric with respect to a plane orthogonal to the axis of
revolution{\rm .}
 \vglue1pt
{\rm (ii)}  If $\th_1\le\pi/6<\th_2\le\pi/3$ then
$\th_2>\displaystyle\frac{\pi}{3}-\th_1${\rm .}
\endproclaim

\figin{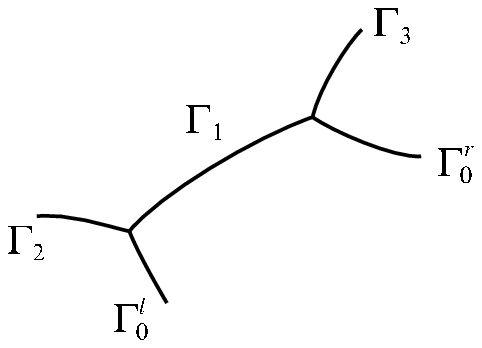}{750}
\vglue2pt
\centerline{Figure 6.}
\vfill

\figin{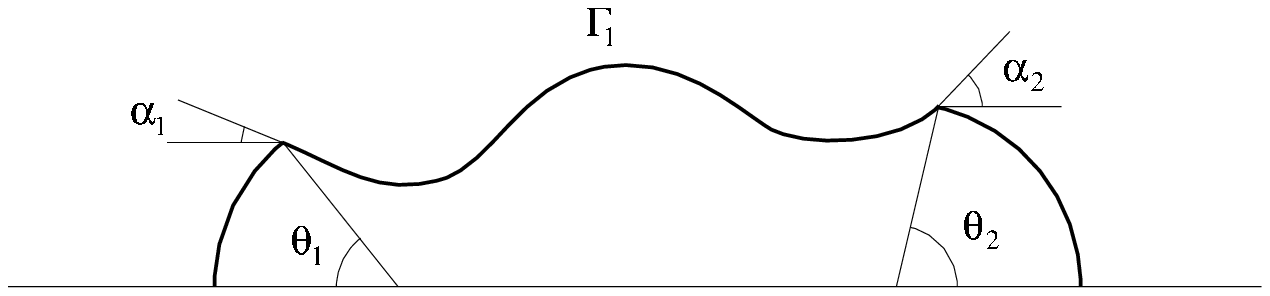}{750}
\vglue2pt
\centerline{Figure 7. The spherical component $\Om$}

\demo{Proof}
Assume that the spherical component $\Om$ is contained in $R_1$.  By
scaling, normalize $H_1=1$.

Let $\Ga_1$ be the generating curve of $M_0$ parametrized from the
left to the right.  As $\alpha_i=\th_i-\pi/6$ we get from
\eqref{eq:fint} that the force $F_0$ of $\Ga_1$ is given by
$g(\th_1)=g(\th_2)$, where
$$
g(\th)=\left(\frac{1}{2}-H_0\right)\,\sin^2\th+\frac{\sqrt{3}}{2}\,\cos\th\,\sin
\th.
$$

Since $H_0=1-H_2<1$ by Proposition~\ref{th:hu}, we get
$g'(\th)\ge-\frac{1}{2}\,\sin 2\th+\frac{\sqrt{3}}{2}\,\cos 2\th>0$,
and so the function $g$ is strictly increasing in $[0,\pi/6]$.  Hence
$\th_1$, $\th_2\le\pi/6$ implies $\th_1=\th_2$.  Moreover the
endpoints of $\Ga_1$ have the same height $y_i$ and the same angle
$\alpha_i$.  This proves (i) by the uniqueness for solutions to
\eqref{eq:generating} with respect to the initial conditions.

To see (ii) let
$$
h(\th)=-\frac{1}{2}\,\sin^2\th+\frac{\sqrt{3}}{2}\,\sin\th\,\cos\th=
g(\th)-(1-H_0)\,\sin^2\th.
$$

It is easily proved that the function $h$ is symmetric with respect to
$\pi/6$ and increasing in $[0,\pi/6]$.  Thus we have
\begin{eqnarray*}
g\biggl(\frac{\pi}{3}-\th_2\biggr)
&=&(1-H_0)\,\sin^2\biggl(\frac{\pi}{3}-\th_2\biggr)
+h\biggl(\frac{\pi}{3}-\th_2\biggr) \\
&<&
(1-H_0)\,\sin^2\th_2+h(\th_2)=g(\th_2)=g(\th_1)
\end{eqnarray*}
and, as $g$ is increasing in $[0,\pi/6]$, we get $(\pi/3)-\th_2<\th_1$, as
we wished.
\enddemo
\pagebreak

\numbereddemo{{R}emark}
If $\Sg$ were $n$-dimensional, then the force of $\Ga_1$ would be
given by $F_0 = g(\th)\,\sin^{n-2}\th$ and so
Lemma~\ref{lem:delaunaygraph} works in arbitrary dimension.
\enddemo
 
\vglue-12pt
\section{Separation and instability}
\setcounter{equation}{0}

Let $\Sigma\subset\R^{n+1}$ a stationary double bubble of revolution
whose axis $L$ is the $x_1$-axis with generating curve
$\Gamma\subset\{(x_1,x_2) |\ x_2\ge 0\}$ consisting of circular arcs
$\overline{\Gamma}_0$ meeting the axis and other arcs
$\overline{\Gamma}_i$ meeting in threes, with interiors $\Gamma_i$
(see Figure~9).  The bubble $\Sigma$ is invariant under the
action of the group $O(n)$ of orthogonal transformations in $\R^{n+1}$
which fix the $x_1$-axis.  We consider the map $f:\Ga - L
\longrightarrow L\cup \{\infty\}$ which maps each $p\in \Ga - L$ to
the point $L(p)\cap L$, where $L(p)$ denotes the normal line to
$\Gamma$ at $p$.  If $L(p)$ does not meet $L$, we define the image of
$p$ as $f(p)=\infty$.  Note that $f$ is multivalued at the endpoints
of the arcs $\Ga_i$, where three of them meet.  We will use the
notation $iA$ and $iB$ for the image under $f$ of the endpoints of
$\overline{\Ga}_i$.

\numbereddemo{{R}emark}
Using \eqref{eq:generating} and \eqref{eq:fint}, we find that the
derivative of $f$ with respect to arc length is given by
$f'=\frac{nF}{y^{n-1}\cos^2\alpha}$.  In particular, $f$ is increasing
as we move to the right along an unduloid or the concave up portion of
a nodoid, decreasing as we move to the right along the concave down
portion of a nodoid, and of course constant on spheres and vertical
hyperplanes.  Hence $f$ is locally injective on any Delaunay curve
with nonzero force.
\enddemo

\proclaim{Proposition}
\label{prop:unstable}
Consider a stable double bubble of revolution $\Sigma\subset\R^{n+1}${\rm ,}
$n\ge 2${\rm ,} with axis $L${\rm .}  Assume there is a finite number of points
$\{p_1,\ldots, p_k\}$ in $\bigcup_i\Ga_i$ with $x = f(p_1) = \cdots =
f(p_k)$ which separates $\Gamma${\rm .}  Assume further that $\{p_1,\ldots,
p_k\}$ is a minimal set with this property{\rm .}

Then every connected component of $\Sg_0${\rm ,} $\Sigma_1$ and $\Sigma_2${\rm ,}
which contains one of the points $p_i$ is part of a sphere centered at
$x$ $($if $x\in L)$ or a part of a hyperplane orthogonal to $L$ $($if
$x=\infty)${\rm .}
\endproclaim

\demo{Proof}
First suppose that $x\in L$ and take, after translation, $x=0$.  The\break
$1$-parameter group of rotations $$\varphi_\theta (x_1,\ldots,x_{n+1})
= (\cos\theta\,x_1 + \sin\theta\,x_2,-\sin\theta\,x_1
+\cos\theta\,x_2,x_3,\ldots,x_{n+1}),$$ $\theta\in \R$, determines a
Jacobi function on the bubble, $u:\Sigma\rightarrow\R$, given by
$$
u(p) = \left< \rest{\frac{d}{d\theta}}{\theta =0}\varphi_\theta(p) ,
N(p)\right> =
-\det (p,N(p),e_3,\ldots,e_{n+1}),
$$
where $N(p)$ is the unit normal vector of $\Sigma$ at $p$,
$\{e_1,\ldots,e_{n+1}\}$ is the standard orthonormal basis of
$\R^{n+1}$ and $\det$ denotes the Euclidean volume element.  We define
here $M_0=\Sigma\cap\{x_2=0\}$.  By the symmetry of $\Sigma$, if $p\in
M_0$, then the vector $N(p)$ also lies in the hyperplane $x_2=0$ and
therefore $u$ vanishes on $M_0$.

On the other hand, if we take $p$ in $f^{-1}\{ 0\}$, then the vectors
$N(p)$ and $p$ are collinear.  Using again the invariance of $\Sigma$
with respect to $O(n)$, we get that $u$ vanishes on the orbit $M (p)$
of $p$ under the action of $O(n)$ (note that $M (p)$ is a hypersurface
of $\Sigma$).

As the points $p_1,\ldots,p_k$ separate the curve $\Gamma$, the set
$M(p_1)\cup\ldots\cup M(p_k)\cup M_0$ is a hypersurface of the bubble
contained in $u^{-1}\{0\}$ which separates $\Sigma$ in at least four
connected components.  In fact, as the set $\{p_1,\ldots,p_k\}$ is
minimal among the subsets of $f^{-1} \{0\}$ satisfying the separation
property, it follows that $\Sigma - [M(p_1)\cup\ldots,M(p_k)\cup M_0]$
has exactly four components $\Lambda_1,\ldots,\Lambda_4$ and that each
one of the sets $M(p_1), \ldots, M(p_k), M_0$ meets the boundary of
each one of these four components.

Consider the functions $v^{(i)}$, $i=1,\ldots,4$, on $\Sigma$ given by
$$
v^{(i)} =
\left\{ \begin{array}{ll}
u, & \hbox{ on } \Lambda_i, \\
0, & \hbox{ on } \Sigma -\Lambda_i.
\end{array}\right.
$$
Then $v^{(i)}$ are admissible and we can find scalars $a_1$, $a_2$,
$a_3$, not all equal to zero, such that $v=\sum_{i=1}^3 a_i v^{(i)}$
verifies the mean value conditions \eqref{mvcond}, so that $v\in {\cal
F}(\Sigma)$.  By Lemma~\ref{Qww},
$$
Q(v,v) = \sum_{i=1}^{3} a_i^2\,Q(v^{(i)},v^{(i)}) =0.
$$
Since $u$ is a Jacobi function,
\begin{equation}
\label{eqn:deltav}
\Delta v+|\sigma|^2v=0
\end{equation}
on $\Sigma\setminus [M(p_1)\cup\cdots\cup M(p_k)\cup M_0]$.  By our
stability hypothesis and Lemma~\ref{Qe0}, equation \eqref{eqn:deltav}
holds on all of $\Sigma$.

Fix $i$ and let $S$ be the connected component of a smooth piece of
$\Sigma$ which contains the point $p_i$.  As $p_i$ lies in the
interior of $S$, the four domains $\Lambda_i$ meet the interior of
$S$.  As $v$ vanishes on $S\cap\Lambda_4$, from the unique
continuation property, we conclude that $v=0$ on $S$.  Hence $u=0$ on
$S\cap\Lambda_j$, for any $j\in\{1, 2, 3\}$ such that $a_j\neq 0$.  As
such $j$ exists we conclude that $u=0$ on $S$ again from the unique
continuation property.  Thus the $1$-parameter group of rotations
$\varphi_\theta$ preserves $S$.  Since $S$ is rotationally symmetric
around the $x_1$-axis, we conclude that this component is a part of a
sphere centered at the origin.

This finishes the proof of the proposition if $x$ is a finite point of
the axis~$L$.

It remains to prove the result when $x=\infty$.  In order to prove it
we repeat the argument by considering, instead of the rotations
$\varphi_\theta$, the $1$-parameter group of translations $T_\theta
(x_1,\ldots,x_{n+1}) = (x_1,x_2 + \theta,\ldots,x_{n+1})$ and its
associated Jacobi function $u(p) = \langle N(p),e_2\rangle$.
\enddemo

\proclaim{{C}orollary}
\label{cor:onecomponent}
There is no stable double bubble of revolution in $\rr^{n+1}$ in which
the graph structure is the one in Figure {\rm 8.}
\endproclaim

\begin{center} \BoxedEPSF{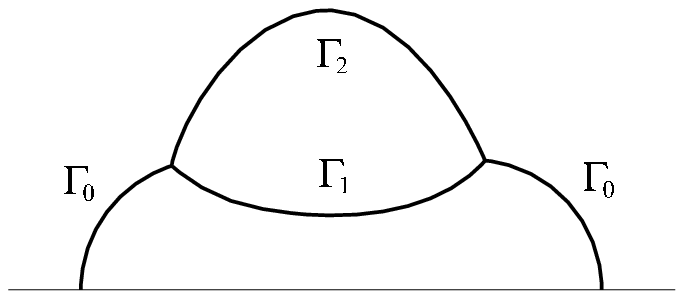 scaled 700} \end{center}
\begin{quote} Figure 8. There is no stable nonstandard double bubble with connected
regions.
\end{quote}

\demo{Proof}
Assume first that the line equidistant from the two vertices
intersects the axis $L$ in a point $p$.  Then $\Ga_1$ and $\Ga_2$ each
has an interior point farthest from or closest to $p$, so that $p\in
f(\Ga_1)\cap f(\Ga_2)$.  By Proposition~\ref{prop:unstable}, $\Ga_1$
and $\Ga_2$ are both spherical (centered on the axis), which is
impossible.

If the equidistant line is horizontal $\Ga_1$, $\Ga_2$ each has an
interior point farthest left or right, so that $\infty\in f(\Ga_1)\cap
f(\Ga_2)$.  By Proposition~\ref{prop:unstable}, $\Ga_1$ and $\Ga_2$
are both vertical, which is impossible.
\enddemo

\numbereddemo{{R}emark}
When $n=2$ and the volumes are equal, Hutchings \cite[Thm.~5.1,
Cor.~4.4]{Hu} showed, as described in our Section~6, that any
nonstandard minimizing bubble satisfies the hypotheses of
Corollary~\ref{cor:onecomponent}.  Therefore in this case the
minimizing solution is the standard bubble.  This fact was first
proved by computer analysis by Hass and Schlafly~\cite{HS2}.
\enddemo

\proclaim{{C}orollary}
\label{cor:stability}
Consider a stable double bubble of revolution in $\rr^{n+1}$ in which
both regions have positive pressure{\rm .}  Assume that one of the regions
$R_2$ is connected{\rm ,} that the other one $R_1$ has two components and
that the graph structure is the one in Figure~{\rm 9.}

Then there is no $x\in L$ such that $f^{-1}(x)-\Ga_0$ contains points
in the interiors of distinct\/ $\Ga_j$ which separate $\Ga${\rm .}
\endproclaim

\begin{center} \BoxedEPSF{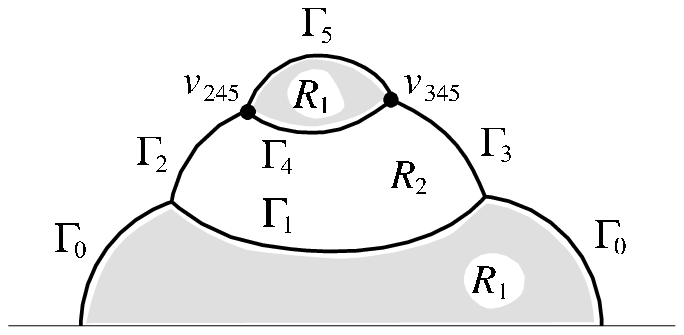 scaled 700} \end{center}

\centerline{Figure 9. A candidate double bubble of three components.}
 
\demo{Proof}
There must be points on $\Gamma_1$, $\Gamma_2$, or $\Gamma_3$.  By
Proposition~\ref{prop:unstable}, one of them is spherical.  By force
balancing~\ref{lem:forcebalancing}, all three are spherical, which is
impossible by Lemma~\ref{lem:gammas}.
\enddemo

We give the following version of Corollary~\ref{cor:stability},
although we will not use it here.

\proclaim{{C}orollary}
\label{cor:new}
For a stable double bubble of revolution{\rm ,} if $f$ is not injective on
the interior of $\Ga_i${\rm ,} then $\Ga_i$ is a circular arc or vertical
line{\rm .}
\endproclaim

\proclaim{Proposition}
\label{prop:separate}
Consider a minimizing nonstandard double bubble in $\rr^{n+1}${\rm ,} which
is necessarily rotationally symmetric around an axis $L${\rm .}

Then there exists no $x\in L$ such that $f^{-1}(x)-\Ga_0$ contains
points in the interiors of distinct\/ $\Ga_j$ which separate $\Ga${\rm .}
\endproclaim

\demo{Proof}
Since every component borders the exterior, a separating set must cut
the outer boundary of some component.  By force
balancing~\ref{lem:forcebalancing}, every piece of the outer boundary
of this component is spherical.  If $f^{-1}(x)$ cuts two pieces of the
outer boundary, then these are pieces of spheres with the same center
and the same mean curvature, and hence the same sphere.  The portion
of the bubble between these two pieces can then be rolled around the
sphere, without changing perimeter or enclosed volume, until it
touches some other part of the bubble, resulting in a bubble which is
not regular, and hence not minimizing.  So it cuts an inner boundary
(part of $\Sg_0$).  By force balancing~\ref{lem:forcebalancing}, each
end of the inner boundary meets two other spheres, which contradicts
Lemma~\ref{lem:gammas}.
\enddemo

\proclaim{Proposition}
\label{prop:twocomponents}
There is no stable double bubble of revolution in $\rr^{n+1}$ in which
both regions have positive pressure{\rm ,} the region of smaller or equal
pressure $R_2$ is connected{\rm ,} the other region $R_1$ has two
components{\rm ,} and the graph structure is the one in Figure~{\rm 9.}
\endproclaim

\demo{Proof}
Suppose there were.  $\Ga_0$ are spherical.  $\Ga_1$ is an unduloid or
catenoid and graph by Lemma~\ref{lem:forcegamma1}.  By force
balancing~\ref{lem:forcebalancing}, $\Ga_2$ and $\Ga_3$ are (convex)
nodoids.  Since the top, third component has larger pressure, $\Ga_4$
must be a (convex) nodoid, catenoid, or vertical line, unless it is
upside down (which cannot occur in the principal cases of
Figure~14).  (Here by ``convex'' we just mean that the tangent
vector rotates monotonically.)

We focus on the third component and its two vertices $v_{245}$ and
$v_{345}$.  For the simplest case when all the curves are graphs as in
Figure 10A, then the images $iA$ and $iB$ under $f$ of the
left and right endpoints of $\Ga_i$ satisfy
\vglue2pt
\centerline{$
4A < 2B < 5A\qquad\hbox{ and}\qquad 5B < 3A < 4B.
$}
\smallbreak\noindent
This remains true as a vertex rotates until one of the three tangent
vectors goes vertical.  (The borderline position with $5A = \infty$
may be considered an extreme position of either case; in the proof we
consider it part of the second case eliminated.)  Rotating $v_{245}$
counterclockwise one notch as in Figure 10B yields instead $5A
< 4A < 2B$.  To avoid giving $\Gamma_4$ or $\Gamma_5$ two vertical
tangents contrary to Corollary~\ref{cor:new}, the two vertices must be
rotated in the same direction, say counterclockwise.  Suppose that
$v_{245}$ is rotated $m_1$ notches counterclockwise and that $v_{345}$
is rotated $m_2$ notches counterclockwise.  Then $m_1 \le 2$, or
$\Gamma_2$ (where $R_2$ is a convex region by positive pressure) could
not meet the circle $\Ga_0$ at $120$ degrees (see Figure 11).
Also $m_2 \le 3$, or $\Gamma_4$ would go vertical twice (see Figure
12), contrary to Corollary~\ref{cor:new}.
\pagebreak

\begin{center} \BoxedEPSF{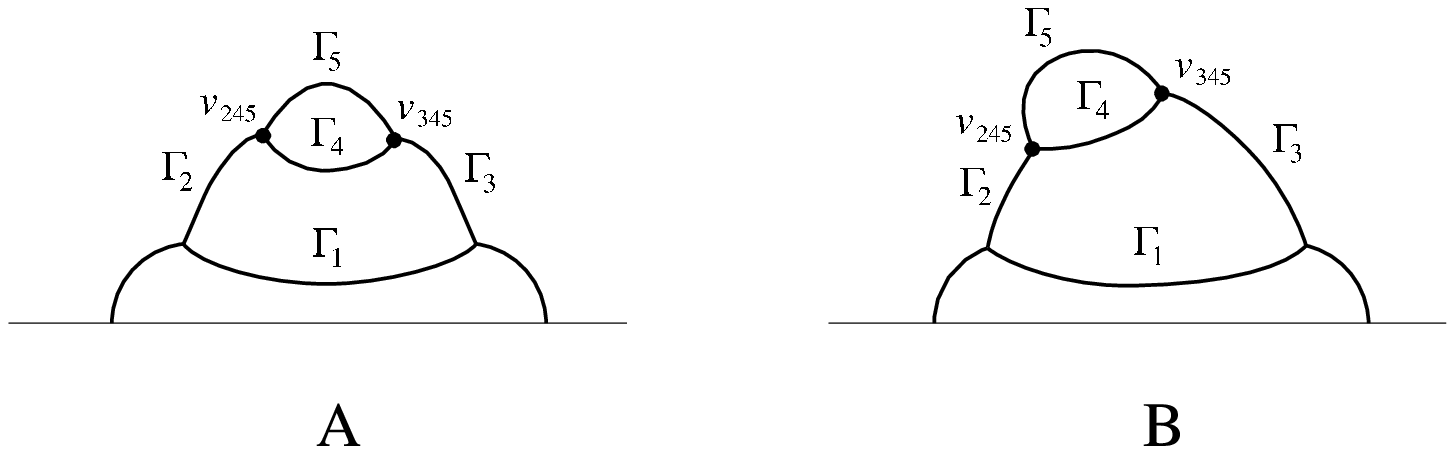 scaled 730} \end{center}
\begin{quote}
Figure 10. Rotating vertex $v_{245}$ one ``notch" (counterclockwise)
means turning one tangent (here $\Ga_5$) past vertical.  (Rotating
another notch would turn $\Ga_2$ past vertical.)
\end{quote}

\begin{center} \BoxedEPSF{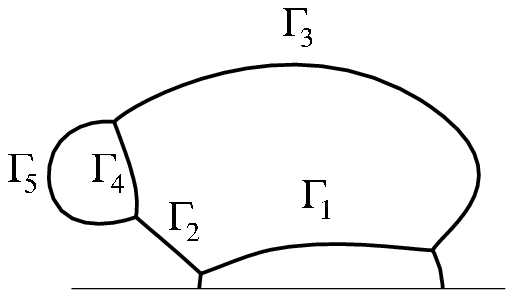 scaled 730} \end{center}
\centerline{Figure 11. If $v_{245}$ turns three notches, $\Ga_2$ cannot be convex.}

\figin{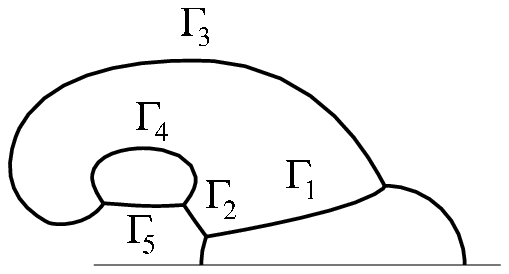}{750}
\centerline{Figure 12. If $v_{345}$ turns four notches, $\Ga_4$ goes vertical twice.}
\vglue4pt

Next consider cases with $m_2 = 2$ or $m_2 = 3$, as in Figure
13.  $\Gamma_3$ is not a graph, $f(\Gamma_3) = [\infty, 3A)
\cup (3B, \infty]$, and by stability Proposition \ref{prop:unstable}
gives that, for $\Gamma_3$, $3A \le 3B$.  We then have $3B<1B$, or
else $\Gamma_3$ would go vertical a second time near $3B$,
contradicting Corollary~\ref{cor:new}.  Therefore $1B$ is contained in
$f(\Gamma_3)$, which contradicts Corollary~\ref{cor:stability} for
$\Gamma_1$ and $\Gamma_3$. 

\begin{center}\BoxedEPSF{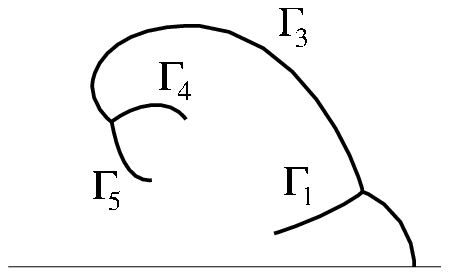 scaled 730} \end{center}
\begin{quote} Figure 13. If $m_2 = 2$ or $m_2 = 3,$ then $\Ga_3$ is not a graph and
$f(\Gamma_3) = [\infty, 3A) \cup (3B, \infty]$.
\end{quote}

Now remain only the cases $0 \le m_1 \le 2$, $0 \le m_2 \le 1$ of
Figure 14.  (It is easy to see that $\Ga_4$ cannot be a
vertical line.)  We claim that
\begin{equation}
\label{eq:(1)}
3A < f(\Gamma_1).
\end{equation}
This is easy if $\Gamma_3$ is a graph, since consideration of $v_{13}$
shows that $3B < 1B$, and then $f(\Gamma_3) < f(\Gamma_1)$ by
Corollary~\ref{cor:stability} for $\Gamma_1$ and $\Gamma_3$.  $v_{13}$
can rotate only clockwise one notch to keep the stem part of a circle
and $\Gamma_1$ a graph.  Now $\Gamma_3$ is not a graph and
$f(\Gamma_3)$ includes $[\infty, 3A)$.  By
Corollary~\ref{cor:stability} for $\Gamma_1$ and $\Gamma_3$, $3A <
f(\Gamma_1)$ as claimed.

\begin{center}\BoxedEPSF{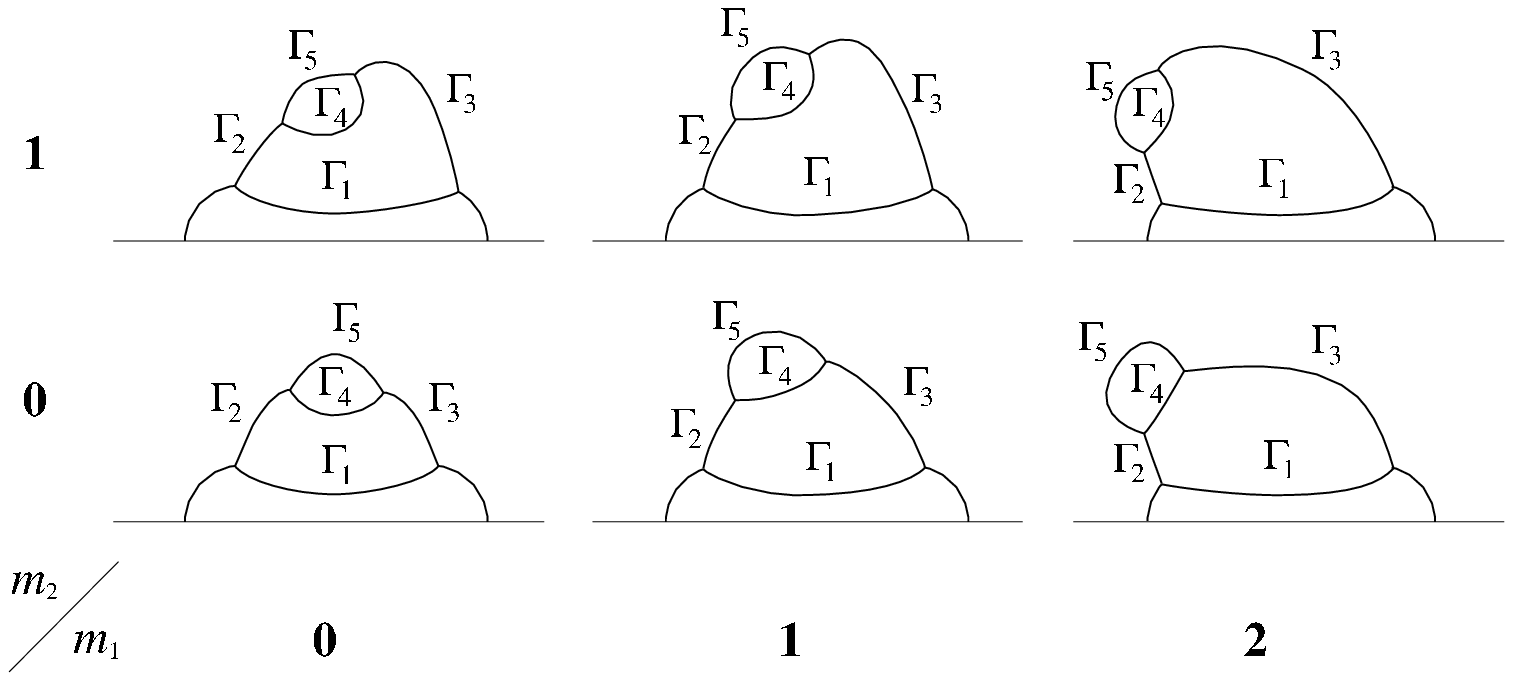 scaled 750} \end{center}
\centerline{Figure 14. The six principal cases to be eliminated}
\vglue4pt
 
For the cases $(0,0)$, $(0,1)$, a similar argument shows that
$f(\Gamma_1)< 2B$.  Consideration of the vertices $v_{245}$ and
$v_{345}$ leads to the conclusion
\begin{equation}
\label{eq:(2)}
5B < 3A < f(\Gamma_1) < 2B < 5A.
\end{equation}
Since the net angle $\theta_5$ through which $\Gamma_5$ turns
satisfies $\theta_5 \le 180$ degrees, obviously $4A < 3A$.  For the
case $(0,1)$, where $f(\Gamma_4)$ contains $(4A, \infty]$, this puts
$3A$ in $f(\Gamma_4) \cap f(\Gamma_5)$, a contradiction of
Corollary~\ref{cor:stability} for $\Gamma_3$, $\Gamma_4$, $\Gamma_5$.
For the case $(0,0)$, consideration of $v_{345}$ shows that $3A < 4B$
and leads to the same contradiction. \pagebreak

Next we consider the cases $(1,0)$, $(2,0)$.  Since $5B < 3A$, $3A$ is
contained in $f(\Gamma_5)$.  Since $3A < 4B$, by
Corollary~\ref{cor:stability} for $\Gamma_3$, $\Gamma_4$, $\Gamma_5$,
we must have $3A \le 4A$.  In particular, $\theta_5 > 180$ degrees.
Of course by Corollary~\ref{cor:new} for $\Gamma_5$, $5A \le 5B$.
Moreover $\Gamma_5$ leaves $v_{345}$ above the horizontal.  Now
Corollary~\ref{cor:final} implies that $3A > 4A$, a contradiction.

Similarly for the final cases $(1,1)$ and $(2,1)$, $3A$ is contained
in $f(\Gamma_5)$.  Since $f(\Gamma_4)$ includes $(4A, \infty]$, by
Corollary~\ref{cor:stability} for $\Gamma_3$, $\Gamma_4$, $\Gamma_5$,
we must have $3A \le 4A$, an immediate contradiction in case $(1,1).$
In particular, $\theta_5 > 180$ degrees, and $5A \le 5B$.  If
$\Gamma_5$ leaves $v_{345}$ at or above the horizontal,
Corollary~\ref{cor:final} yields the contradiction $3A > 4A$.  If on
the other hand $\Gamma_5$ leaves $v_{345}$ below the horizontal, then
the downward normal $n$ to $\Gamma_3$ at $v_{345}$ is counterclockwise
from the downward tangent to $\Gamma_2$ at $v_{12}$ (and hence from
every downward tangent to $\Gamma_2$) and hence counterclockwise from
the downward normal to $\Gamma_1$ at $v_{12}$.  Since $\Gamma_4$ is
convex, $n$ stays to the right of $\Gamma_2$ and $1A < 3A$, a
contradiction of \eqref{eq:(1)}.
\enddemo

\proclaim{Lemma}
\label{lem:trig}
Given points $A$ and $B${\rm ,} consider two points $D${\rm ,} $E$ on the same
side of $AB$ subtending the same angle $\theta$ as in
Figure~{\rm 15.}  Then $\angle CDE = \angle ABC${\rm .}
\endproclaim

\demo{Proof}
Since $\angle BCA = \angle DCE$ and $ACE \sim BCD$, $ABC \sim CDE$.
\enddemo

\begin{center}\BoxedEPSF{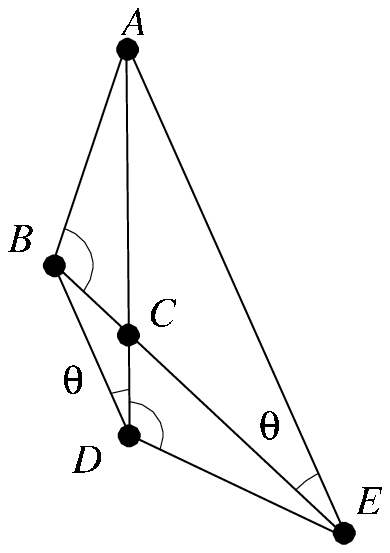 scaled 800} \end{center}
\centerline{Figure 15. $\angle CDE=\angle ABC$}
  \vglue4pt 

\proclaim{{C}orollary}
\label{cor:final}
In cases $(1,0)${\rm ,} $(2,0)${\rm ,} and $(2,1)$ of Figure~{\rm 14,} suppose
that the net angle $\theta_5$ through which $\Gamma_5$ turns exceeds
$180$ degrees{\rm ,} that $\Gamma_5$ leaves $v_{345}$ at or above the
horizontal{\rm ,} and that $5A \le 5B${\rm .}  Then $4A < 3A${\rm .}
\endproclaim

\demo{Proof}
Let $\theta = \theta_5 - 180 > 0$.  Apply Lemma~\ref{lem:trig} with $A
= v_{345}$, $B =v_{245}$, $AD \perp \Gamma_5$ and $AE \perp \Gamma_3$
at $v_{345}$, $BD \perp \Gamma_5$ and $BE \perp \Gamma_4$ at
$v_{245}$; then $\angle ADB = \angle AEB = \theta$.  Since $\Gamma_4$
is strictly convex (it cannot be a vertical line because $5A \le 5B$),
$\angle ABC > 90$.  By Lemma~\ref{lem:trig}, $\angle CDE = \angle ABC
> 90$.  Since by hypothesis $\Gamma_5$ leaves $v_{345}$ at or above
the horizontal, $DE$ heads downward.  (In these cases,
Figure~15, in which $AD$ is vertical, is rotated clockwise
by an amount less than $90$ degrees, strictly less by
Corollary~\ref{cor:new} since $\theta_5 > 180$.)  Since $5A \le 5B$,
$D$ lies on or below the horizontal axis.  Hence $E$ lies below the
horizontal axis and $4A < 3A$.
\enddemo 

\section{Estimates on the number of components}
\label{sec:comp}
 
In this section we prove that for a minimizing double bubble in
$\rr^3$, the larger region is connected, and the smaller region has at
most two components.

We begin by recalling a version of the Hutchings Basic Estimate.  Let
$A(v)$ denote the volume of a sphere in $\rr^3$ enclosing volume $v$,
and let $A(v_1,v_2)$ denote the area of the standard double bubble
enclosing volumes $v_1$ and $v_2$ in $\rr^3$.  We then have:

\proclaimtitle{{\cite[Thm.~4.2]{Hu}}}
\proclaim{Proposition}
\label{prop:last}
In an area\/{\rm -}\/minimizing double bubble enclosing volumes $v_1$ and $v_2$
in $\rr^3${\rm ,} suppose that $R_2$ contains a component of volume $\lambda
v_2${\rm .}  Then
$$
2A(v_1,v_2) \ge \lambda^{-1/3}A(v_2)+A(v_1)+A(v_1+v_2).
$$
\endproclaim

This inequality places a lower bound on $\lambda$ in terms of $v_1$
and $v_2$.  Hutchings \cite[4.4, 4.5]{Hu} calculated that when
$v_1=v_2$, the lower bound is greater than $1/2$, so both regions are
connected; and when $v_1>>v_2$, the lower bound approaches $2/5$, so
the smaller region has at most $2$ components.  More recently,
Heilmann et al.  \cite[Fig.\  8]{HLRS} (\cite[Fig.\ 14.11.1]{M2}) used a
computer to plot the lower bound on $\lambda$ as a function of the
ratio $v_2/v_1$, and found that it is apparently increasing.  This
would imply that the larger region is always connected and the smaller
region always has at most two components.  However this observation
has not been rigorously proved, because the function $A(v_1,v_2)$ is
difficult to work with.  Thus we will use different methods to prove
the above connectedness results.

\proclaim{Proposition}
\label{prop:largeconnected}
In a minimizing double bubble in $\rr^3${\rm ,} the region with larger or
equal volume is connected{\rm .}
\endproclaim

\demo{Proof}
By rescaling, we may assume that the two volumes are $1-w$ and~$w$.
Hutchings \cite[Thm.  3.5]{Hu} showed that if $w<1/3$, then the larger
region is connected.  (This is true in higher dimensions as well.)

For $w\ge 1/3$, to prove that the region of volume $1-w$ is connected,
it is enough to show that the inequality of
Proposition~\ref{prop:last} fails for $\lambda=1/2$; i.e.
$$
2A(w,1-w) < 2^{1/3}A(1-w)+A(w)+A(1).
$$
We observe that by Lemma~\ref{lem:1starea}, $dA(w,1-w)/dw>0$ for
$w<1/2$, because we can continuously deform one standard double bubble
to another, and the smaller region has larger pressure.  Thus
$A(w,1-w) \le A(1/2,1/2)$.  It is straightforward to compute that
$A(1/2,1/2)=2^{-4/3}3A(1)$, and $A(v)=v^{2/3}A(1)$.  Thus it will
suffice to show that
$$
2^{-1/3}3 < 2^{1/3}(1-w)^{2/3}+w^{2/3}+1.
$$
Since this holds at $w=0.10$ and $w=0.63$ and the right-hand side is
concave, it holds for $0.10\le w\le 0.63$.  In particular, it holds
for $1/3\le w \le 1/2$, as desired.
\enddemo

\numbereddemo{{R}emark}
An alternative proof is given by Heilmann et al.
\cite[Prop.~2.5]{HLRS}.  Actually as in the proof above, any region
with at least $37\%$ of the total volume is connected.
\enddemo

\proclaim{Lemma}
\label{lemma:ros}
In a minimizing double bubble in \/ $\rr^{n+1}$ enclosing two unequal
volumes{\rm ,} the smaller region has larger pressure{\rm .}
\endproclaim

\demo{Proof}
Consider the function $A(v,1-v)$ giving the least area enclosing and
separating regions of volume $v$, $1-v$.  By Hutchings
\cite[Thm.~3.2]{Hu}, $A$ is strictly concave and of course symmetric
about $v=1/2$.  Moving the separating surface (of mean curvature
$H_0$) of an area-minimizing bubble we have $dA/dv=nH_0$, and the left
and right derivatives of $A$ must satisfy
$$
A_R'\le nH_0 \le A_L'.
$$
Consequently $H_0$ is positive for $v<1/2$ and negative for $v>1/2$.
In other words, the smaller region has larger pressure.
\enddemo

The following Proposition~\ref{prop:small2comp} shows that small
region with three components is unstable, by using a volume-preserving
linear combination of (non\-volume-preserving) variations of each of
the components.

\proclaim{Proposition}
\label{prop:small2comp}
In a minimizing double bubble in $\rr^3${\rm ,} the region with smaller or
equal volume has at most two connected components{\rm .}
\endproclaim

\demo{{P}roof}
Assume that the volume of $R_1$ is less than or equal to the volume of
$R_2$.  By Proposition~\ref{th:hu} we obtain $H_1$, $H_2>0$.  By
Lemma~\ref{lemma:ros} we get $H_0\ge 0$.

Recall that $\kappa_i=\sg_i(\nu_i,\nu_i)$, and let $c_i=\sg_i(T,T)$,
where $T$ is the unit tangent vector to the singular curve $C$.  So
$2H_i=\kappa_i+c_i$.

We consider an admissible function $u$ invariant by the one-parameter
group of rotations of $\Sg$.  The functions $u_i$ are locally constant
over $C$.  If we apply \eqref{eq:defq} to $u$, adding and subtracting
$c_iu_j^2$ in the boundary term, \pagebreak we see that the second variation form
satisfies
\begin{eqnarray*}
Q(u,u)&=&\sum_i \int_{\Sg_i} \{|\nabla u_i|^2-|\sg_i|^2 u_i^2\} \\
&&-\ \frac{2}{\sqrt{3}} \int_C
(H_0-H_2)\,u_1^2+(-H_1-H_0)\,u_2^2+(H_1+H_2)\,u_0^2 \\
&&+\ \frac{1}{\sqrt{3}} \int_C
(c_0-c_2)\,u_1^2+(-c_1-c_0)\,u_2^2+(c_1+c_2)\,u_0^2.
\end{eqnarray*}
Taking the scalar product with $D_T T$ in the formulae \eqref{eq:nui}
we have
$$
\kappa_g^1=\frac{1}{\sqrt{3}}\,(c_0-c_2), \qquad
\kappa_g^2=\frac{1}{\sqrt{3}}\,(-c_1-c_0), \qquad
\kappa_g^0=\frac{1}{\sqrt{3}}\,(c_1+c_2),
$$
where $\kappa_g^i$ stands for the geodesic curvature of $C$ inside
$\Sg_i$ (with respect to the conormal $\nu_i$).  So we have
\begin{eqnarray}
\label{eq:quu}\qquad
Q(u,u)&=&\sum_i \int_{\Sg_i} \left\{|\nabla u_i|^2- |\sg_i|^2 u_i^2\right\}
\\
&&-\ \frac{2}{\sqrt{3}} \int_C
(H_0-H_2)\,u_1^2+(-H_1-H_0)\,u_2^2+(H_1+H_2)\,u_0^2 \nonumber\\
 &&+\ \int_C
\kappa_g^1\,u_1^2
+\kappa_g^2\,u_2^2
+\kappa_g^0\,u_0^2.\nonumber
\end{eqnarray}

Consider a connected component $\Om$ of $R_1$.  Let
$M_i=\Sg_i\cap\ptl\Om$, and let $C^*=C\cap\ptl\Om$.  We want to find
an admissible function $u$ such that $Q(u,u)<0$ with support inside
$\ptl\Om$.  Then if $R_1$ had three connected components, some nonzero
linear combination would preserve the volumes and yield a
contradiction.

We define the function
$$
v=\left\{ \begin{array}{ll}
1, & \hbox{ on}\ M_0\cup M_1, \\
0, & \hbox{ elsewhere in}\ \Sg.
\end{array}\right.
$$
Then \eqref{eq:quu} gives
\begin{equation}
\label{eq:q11a}
Q(v,v)=-\sum_{i=0, 1}\int_{M_i} |\sg_i|^2-\frac{2}{\sqrt{3}}\int_{C^*}
(H_0+H_1)+\int_{C^*} (\kappa_g^0+\kappa_g^1).\hskip.25in
\end{equation}

Since $|\sg_i|^2=4H_i^2-2K_i$, from \eqref{eq:q11a} and Gauss-Bonnet
$\int_{M_i}K_i=2\pi\chi(M_i)-\int_{\ptl M_i}\kappa_g^i$ we obtain
\begin{equation}
\label{eq:q11b}
Q(v,v)=\sum_{i=0,1} \biggl\{4\pi\chi(M_i)-\int_{M_i}
4H_i^2\biggr\}-\frac{2}{\sqrt{3}}\int_{C^*} (H_0+H_1)-\int_{C^*}
(\kappa_g^0+\kappa_g^1).
\end{equation}

Assume first that $\Om$ is a torus component, so that its boundary is
a union of annuli.  Adding \eqref{eq:q11a} and \eqref{eq:q11b} and
taking into account that \pagebreak $\chi(M_i)=0$, we eliminate the geodesic
curvature to obtain
$$
2 Q(v,v)=-\sum_{i=0, 1}\int_{M_i}
\bigl\{|\sg_i|^2+4H_i^2\bigr\}-\frac{4}{\sqrt{3}}\int_{C^*} (H_0+H_1)<0,
$$
as desired.

We now assume that $\Om$ is the spherical component, so that $M_1$ is
the union of two spherical caps $D_1$, $D_2$ and an annulus $M_0$, as
in Figure~7.  As $M_0$ is a graph by
Lemma~\ref{lem:forcegamma1} we conclude $0<\th_i\le\frac{2\pi}{3}$,
where $\th_i$ is the angle determined by $D_i$.  By scaling we may
assume that the spherical caps have mean curvature $H_1=1$.  Using
Gauss-Bonnet we get that
\begin{equation}
\label{eq:gb}
A(M_1)=\int_{M_1} K_1=4\pi-\int_{C^*}\kappa_g^1.
\end{equation}
Since $\nu_0=(-1/2)\,\nu_1+(\sqrt{3}/2)\,N_1$, taking the scalar
product with $D_T T$ we have
\begin{equation}
\label{eq:kappa}
\kappa_g^0=-\frac{1}{2} \kappa_g^1+\frac{\sqrt{3}}{2}.
\end{equation}

From \eqref{eq:q11b}, \eqref{eq:gb} and \eqref{eq:kappa} we obtain,
taking into account that $\chi(M_1)=2$ and discarding the summands
containing $H_0$,
$$
Q(v,v)\le 6\pi+\frac{7}{2}\,A(M_1)-\frac{7}{2\sqrt{3}}\,L(C^*).
$$
As $A(D_i)=2\pi\,(1-\cos\th_i)$ and $L(\ptl D_i)=2\pi\sin\th_i$, we
have
\begin{equation}
\label{eq:qvv}
Q(v,v)\le 2\pi\left\{-4+h(\th_1)+h(\th_2)\right\},
\end{equation}
where $h(\th)=\frac{7}{2}\,
\left(\cos\th-\frac{1}{\sqrt{3}}\,\sin\th\right)$, which is decreasing
in $[0,2\pi/3]$.  Thus if $\th_1$ or $\th_2$ is greater than or equal
to $\pi/2$, we have
$$
Q(v,v)\le 2\pi\{h(0)+h(\pi/2)\}<0.
$$

Assume now that both $\th_1$, $\th_2<\pi/2$, and consider the function
$$
w=\left\{ \begin{array}{ll}
\displaystyle\frac{\cos\th}{\cos\th_i}, & \hbox{ in}\ D_i, \\
1, & \hbox{ in}\ M_0, \\
0, & \hbox{ elsewhere in}\ \Sg.
\end{array}\right.
$$
As $v$ and $w$ differ only on $M_1$ we obtain from \eqref{eq:defq}
\begin{equation}
\label{eq:qww}
Q(w,w)=\int_{M_1}\left(|\nabla w|^2-2 w^2\right)+2\int_{M_1} 1+Q(v,v).
\end{equation}
By direct computation we get
$$
\int_{M_1}\left(|\nabla w|^2-2
w^2\right)=-2\pi\sum_{i=1,2}\frac{\sin^2\th_i}{\cos\th_i},  \qquad
\int_{M_1} 1=2\pi\,\sum_{i=1,2} (1-\cos\th_i),
$$
which combined with \eqref{eq:qvv} and \eqref{eq:qww} yield
\begin{equation}
\label{eq:qww2}
Q(w,w)\le 2\pi\left\{g(\th_1)+g(\th_2)\right\},
\end{equation}
where $g$ is given by
$$
g(\th)=\frac{3}{2}\cos\th-\frac{7}{2\sqrt{3}}
\sin\th-\frac{\sin^2\th}{\cos\th}.
$$

The function $g$ is strictly decreasing in $[0,\pi/2]$ since it is the
sum of three decreasing functions.  As $g(0)=\frac{3}{2}>0$,
$g(\pi/6)=0$, and $g(\pi/3)=-5/2$ we conclude
$$
Q(w,w)<0\quad \hbox{ if either both}\quad \th_1,
\th_2>\frac{\pi}{6}\quad\hbox{ or some}\quad  \th_i\ge\frac{\pi}{3}.
$$

We finally consider the remaining cases in which at least one of the
angles $\th_i$ is less than or equal to $\pi/6$ and both are less than
$\pi/3$.

\demo{Case  1}
$\th_1\le\frac{\pi}{6}<\th_2<\frac{\pi}{3}$.
\enddemo

Observe that $g$ is concave in $[0,\pi/3]$ since
$$
g''(\th)=-\frac{7}{2}\,\biggl(\cos\th-\frac{1}{\sqrt{3}}\,\sin\th\biggr)-\biggl(
3\,\frac{\sin^2\th}{\cos\th}+2\,\frac{\sin^3\th}{\cos^2\th}\biggr)<0.
$$
By Lemma~\ref{lem:delaunaygraph} we know that $\frac{\pi}{3}-\th_1<\th_2$.
As $g$ is decreasing and concave
$$
\frac{1}{2\pi}\,Q(w,w)\le
g(\th_1)+g(\th_2)<g(\th_1)+g\biggl(\frac{\pi}{3}-\th_1\biggr)\le
2\,g\left(\frac{\pi}{6}\right)=0.
$$

\demo{Case  2}
$\th_1$, $\th_2\le\frac{\pi}{6}$.
\enddemo

By Lemma~\ref{lem:delaunaygraph} $M_0$ is symmetric with respect to a
plane orthogonal to the line of symmetry.  So if
$\kappa_1=\sg_1(\nu_1,\nu_1)\ge 0$ we get from \eqref{eq:defq}
$$
Q(v,v)=-\int_{M_0\cup M_1} |\sg|^2-\int_{C^*}(\kappa_1+\kappa_0)<0.
$$

If $\kappa_1=\sg_1(\nu_1,\nu_1)<0$ then the Gauss curvature of $M_0$
along $C$ is negative.  By Lemma~\ref{lem:forcegamma1} $M_0$ is an
unduloid or a catenoid.  As $\th_i\le\pi/6$ the vectors $\nu_1$, which
are tangent to the generating curve $\Ga_1$ of $M_0$ at their
endpoints, are either horizontal or upper pointing.  Therefore $M_0$
contains a nodal region of the Gauss curvature in its interior, which
implies that $M_0$ is unstable \cite[Thm.~3]{RR},
\cite[Prop.~4.1]{PR}.

So for any component of $R_1$ we have an admissible function $u$
such that $Q(u,u)<0$ with support inside the boundary of the
component.  If we had three connected components in $R_1$ then we
could get an admissible function satisfying the mean value zero
property \eqref{mvcond}, which gives instability of the considered
double bubble, a contradiction.\hfill\qed
\enddemo
\pagebreak

{\it Remark} 6.6.
\begin{itemize}
\item[(a)] An alternative, computational proof of
Proposition~\ref{prop:small2comp}, using Proposition~\ref{prop:last},
is outlined by Heilmann et al.  \cite[Prop.~4.6]{HLRS} (see
\cite[14.11, 14.13]{M2}).  \item[(b)] Reichardt et al.  \cite{RHLS}
extended the arguments of Section 5 to prove the double bubble
conjecture assuming only that one region is connected, thus providing
an alternative to proving Proposition~\ref{prop:small2comp}.
\end{itemize}

\section{Proof of the double bubble conjecture}
 
\proclaim{Theorem}
\label{th:main}
The standard double bubble in $\rr^3$ is the unique area\/{\rm -}\/minimizing
double bubble for prescribed volumes{\rm .}
\endproclaim

\demo{Proof}
Let $\Sg$ be an area-minimizing double bubble.  By
Propositions~\ref{prop:largeconnected} and \ref{prop:small2comp},
either both regions are connected, or the region of larger volume and
smaller pressure is connected and the one of smaller volume and larger
pressure has two components.  By Proposition~\ref{th:hu}, $\Sg$ is
either the standard double bubble or a bubble like the ones in
Figures~8 or  9.  As $\Sg$ is stable, by
Corollary~\ref{cor:onecomponent} and
Proposition~\ref{prop:twocomponents} we conclude that it must be the
standard double bubble.
\enddemo
  
7.1. {\it Immiscible fluid clusters}.  The methods of this paper extend to double clusters in which the three interfaces carry
different costs, so-called immiscible fluid clusters (see M2, Chapt.\ 16]).  We assume that the costs $a_{01},a_{02},a_{12}$
satisfy strict triangle inequalities, such as
$$
\varepsilon_{02} = a_{01}+a_{12}-a_{02}>0.
$$

\proclaim{Theorem}  For nearly unit costs{\rm ,} if the smaller region has at least $37\%$ of the volume{\rm ,} then the standard
double cluster minimizes energy{\rm .}
\endproclaim

{\it Proof sketch}.  Proposition 6.1 has the following generalization to least energy:
\begin{equation}
2E(v_1,v_2)\geq \lambda^{-1/3} \varepsilon_{01} A(v_2)+\varepsilon_{02}A(v_1)+\varepsilon_{12}A(v_1+v_2).
\end{equation}
When the costs $a_{ij}$ and hence the $\varepsilon_{ij}$ are all $1$, (7.1) reduces to Proposition~6.1.  When they are near $1$, the
proof of Proposition 6.2 still shows that both regions are connected.  Now the simple plane geometry of Corollary~5.3 shows that
the cluster must be standard.\hfill\qed
\pagebreak

{\it Remark} 7.3.  For general costs, even for equal volumes, it remains an open question whether the standard double immiscible
fluid cluster minimizes energy.  The above proof applies whenever we know that both regions are connected.  (The more
complicated plane geometry of Proposition 5.8 (or [RHLS]), for the case when the larger region is connected but the smaller region
has two (or more) components, does not immediately generalize, because the generating curves no longer meet at 120 degrees.) 
Unfortunately, for general costs, even for equal volumes, (7.1) does not imply both regions connected.

We thank undergraduate Ken Dennison for raising this question.
\vglue-9pt

\AuthorRefNames  [HMRR]

 \end{document}